\documentclass[a4paper]{article}
\usepackage{theorem,amsmath,amssymb,amscd}
\usepackage[small]{titlesec}

\theoremstyle{change} \allowdisplaybreaks 

\newcommand{\R}{\mathbb{R}}
\newcommand{\C}{\mathbb{C}}
\newcommand{\Z}{\mathbb{Z}}
\newcommand{\Q}{\mathbb{Q}}
\newcommand{\A}{\mathbb{A}}

\newcommand{\SH}{\mathfrak h}
\newcommand{\HH}{\mathbb{H}}

\newcommand{\GL}{\mathrm{GL}}
\newcommand{\GU}{\mathrm{GU}}

\newcommand{\SL}{\mathrm{SL}}

\newcommand{\GSp}{\mathrm{GSp}}
\newcommand{\SSp}{\mathrm{Sp}}
\newcommand{\OF}{\mathfrak{o}}
\newcommand{\p}{\mathfrak{p}}
\renewcommand{\P}{\mathfrak{P}}
\newcommand{\I}{{\rm I}}

\newcommand{\tr}{{\rm tr}}

\newcommand{\mat}[4]{{\setlength{\arraycolsep}{0.5mm}\left[
\begin{array}{cc}#1&#2\\#3&#4\end{array}\right]}}
\newcommand{\qed}{\hspace*{\fill}\rule{1ex}{1ex}}
\def\qdots{\mathinner{\mkern1mu\raise0pt\vbox{\kern7pt\hbox{.}}\mkern2mu
\raise3.4pt\hbox{.}\mkern2mu\raise7pt\hbox{.}\mkern1mu}}

\newtheorem{thm}{Theorem.}[section]
\newtheorem{theorem}{Theorem.}[section]

\newtheorem{lemma}[thm]{Lemma.}

\newtheorem{corollary}[thm]{Corollary.}

\newtheorem{proposition}[thm]{Proposition.}
\newtheorem{rem}[thm]{Remark.}

\begin{document}

\begin{center}
{\Large Steinberg representation of $\GSp_4$: Bessel models and integral representation of $L$-functions}

\vspace{2ex} Ameya Pitale\footnote{American Institute of Mathematics, Palo Alto, CA 94306, {\tt pitale@aimath.org}}

\vspace{3ex}
\begin{minipage}{100ex}\small
{\sc Abstract.} We obtain explicit formulas for the test vector in the Bessel model and derive the criteria for existence and uniqueness for Bessel models for the unramified, quadratic twists of the Steinberg representation $\pi$ of $\GSp_4(F)$, where $F$ is a non-archimedean local field of characteristic zero. We also give precise criteria for the Iwahori spherical vector in $\pi$ to be a test vector. We apply the formulas for the test vector to obtain an integral representation of the local $L$-function of $\pi$ twisted by any irreducible, admissible representation of $\GL_2(F)$. Together with results in \cite{Fu} and \cite{PS2}, we derive an integral representation for the global $L$-function of the irreducible, cuspidal automorphic representation of $\GSp_4(\A)$ obtained from a Siegel cuspidal Hecke newform, with respect to a Borel congruence subgroup of square-free level, twisted by any irreducible, cuspidal, automorphic representation of $\GL_2(\A)$. A special value result for this $L$-function in the spirit of Deligne's conjecture is obtained. \end{minipage}

\end{center}
\section{Introduction}
It is known that the representation of the symplectic group obtained from a Siegel modular form is non-generic, which means that it does not have a Whittaker model. Consequently, one cannot use the techniques or results for generic representations in this case. In such a situation one introduces the notion of a generalized Whittaker model, now called the Bessel model. These Bessel models have been used to obtain integral representations of $L$-functions. It is known that an automorphic representation of $\GSp_4(\A)$, where $\A$ is the ring of adeles of a number field, obtained from a Siegel modular form always has some global Bessel model. For the purposes of local calculations it is often very important to know the precise criteria for existence of local Bessel models and explicit formulas. In this paper, we wish to investigate Bessel models for unramified, quadratic twists of the Steinberg representation $\pi$ of $\GSp_4(F)$, where $F$ is any non-archimedean local field of characteristic zero.

Let us first briefly explain what a Bessel model is. Detailed definitions will be given in Sect.\ \ref{bessel-model-section}. Let $F$ be a non-archimedean  field of characteristic zero. Let $U(F)$ be the unipotent radical of the Siegel parabolic subgroup of $\GSp_4(F)$ and $\theta$ be any non-degenerate character of $U(F)$. The group $GL_2(F)$, embedded in the Levi subgroup of the Siegel parabolic subgroup, acts on $U(F)$ by conjugation and hence, on characters of $U(F)$. Let $T(F) = {\rm Stab}_{\GL_2(F)}(\theta)$. Then $T(F)$ is isomorphic to the units of a quadratic extension $L$ of $F$. The group $R(F) = T(F)U(F)$ is called the Bessel subgroup of $\GSp_4(F)$ (depending on $\theta$). Let $\Lambda$ be any character of $T(F)$ and denote by $\Lambda \otimes \theta$ the character obtained on $R(F)$. Let $(\pi, V)$ be any irreducible, admissible representation of $\GSp_4(F)$. A linear functional $\beta : V \rightarrow\C$, satisfying $\beta(\pi(r)v) = (\Lambda \otimes \theta)(r) \beta(v)$ for any $r \in R(F), v \in V$, is called a $(\Lambda, \theta)$-Bessel functional for $\pi$. We say that $\pi$ has a $(\Lambda, \theta)$-Bessel model if $\pi$ is isomorphic to a subspace of smooth functions $B : \GSp_4(F) \rightarrow \C$, such that $B(rh) = (\Lambda \otimes \theta)(r) B(h)$ for all $r \in R(F), h \in \GSp_4(F)$. The existence of a non-trivial Bessel functional is equivalent to the existence of a Bessel model for a representation. If $\pi$ has a non-trivial $(\Lambda, \theta)$-Bessel functional $\beta$, then a vector $v \in V$ such that $\beta(v) \neq 0$ is called a \emph{test vector} for $\beta$. 

In \cite{PT}, the authors have obtained, for any irreducible, admissible representation $\pi$ of $\GSp_4(F)$,  the criteria to be satisfied by $\Lambda$ for the existence of a $(\Lambda, \theta)$-Bessel functional for $\pi$. Their method involves the use of theta lifts and distributions. The uniqueness of Bessel functionals has been obtained in \cite{N-PS} for many cases, in particular for any $\pi$ with a trivial central character. In \cite{Su}, a test vector is obtained when both the representation $\pi$ and the character $\Lambda$ are unramified. In \cite{Sa}, a test vector is obtained when $F = \Q_p, p$ is odd and inert in the quadratic field extension $L$ corresponding to $T(\Q_p)$, the representation $\pi$ is an unramified, quadratic twist of the Steinberg representation and $\Lambda$ has conductor $1+p\OF_L$. The explicit formulas of the test vector in the above two cases have been used in \cite{Fu} and \cite{Sa} to obtain an integral representation of the $\GSp_4 \times \GL_2$ $L$-function where the $\GL_2$ representation is either unramified or Steinberg.

The main goal of this paper is to obtain explicit formulas for a test vector whenever a Bessel model for the unramified, quadratic twist of the Steinberg representation of $\GSp_4(F)$ exists. In addition to obtaining these formulas, we, in fact, obtain an independent proof of the criteria for existence \emph{and} uniqueness for the Bessel models.  We also give precise conditions on the character $\Lambda$ so that the Iwahori spherical vector in $\pi$ is a test vector in Theorem \ref{existence-theorem}. The methods used here are very different from those in \cite{N-PS} and \cite{PT}. 

When the Iwahori spherical vector is a test vector, we use the explicit formula for the test vector to obtain an integral representation of the local $L$-function $L(s, \pi \times \tau)$ of the Steinberg representation $\pi$ of $\GSp_4(F)$, twisted by \emph{any} irreducible admissible representation $\tau$ of $\GL_2(F)$. This integral involves a function $B$ in the Bessel model of $\pi$ and a Whittaker function $W^\#$ in a certain induced representation of $\GU(2,2)$ related to $\tau$. We wish to remark that in this paper, and other works (\cite{Fu}, \cite{PS1}, \cite{PS2}, \cite{Sa}), the Bessel function $B$ is always chosen to be a ``distinguished" vector (spherical if $\pi$ is unramified and Iwahori spherical if $\pi$ is Steinberg) which has the additional property of being a test vector. With this choice of $B$ we have a systematic way of choosing $W^\#$ (see \cite{PS2}) so that the integral is non-zero and gives an integral representation of the $L$-function. The work so far suggests that to obtain an integral representation for the $L$-function with a general irreducible, admissible representation $\pi$ of $\GSp_4(F)$, we will have to choose $B$ to be both a ``distinguished" vector in the Bessel model of $\pi$ {\it and} a test vector for the Bessel functional. This further highlights the importance of obtaining more information and explicit formulas for test vectors for Bessel models of $\GSp_4(F)$. This is a topic of ongoing work.

The local computation mentioned above, together with the archimedean and $p$-adic calculations in \cite{Fu} and \cite{PS2}, we obtain an integral representation of the global $L$-function $L(s,\pi \times \tau)$ of an irreducible, cuspidal, automorphic representation $\pi$ of $\GSp_4(\A)$, obtained from a Siegel cuspidal newform with respect to the Borel congruence subgroup of square-free level, twisted by \emph{any} irreducible, cuspidal, automorphic representation $\tau$ of $\GL_2(\A)$. When $\tau$ corresponds to an elliptic cusp form in $S_l(N,\chi)$, we obtain algebraicity results for special value of the twisted $L$-function in the spirit of Deligne's conjecture \cite{De1}.

The paper is organized as follows. The first half of the paper deals with the existence and uniqueness of Bessel models. In Sect.\ \ref{pi-steinberg-setup}, we give the basics regarding the non-archimedean setup, Steinberg representation and the Iwahori Hecke algebra. The Steinberg representation is characterized as the only representation with a unique Iwahori spherical vector. We define $B(\Lambda,\theta)^\I$ to be the space of smooth functions on $\GSp_4(F)$ which are right invariant under the Iwahori subgroup $\I$ and transform on the left according to $\Lambda \otimes \theta$. In Sect.\ \ref{double-coset-decomp-section}-\ref{criterion-for-dim-section}, we obtain that $\dim(B(\Lambda,\theta)^\I) \leq 1$ (the uniqueness), the criterion for $\dim(B(\Lambda,\theta)^\I) = 1$ and the explicit formula for the unique (up to scalars) function $B$ in $B(\Lambda,\theta)^\I$. The methods used here are similar to those in \cite{Sa}. In case $\Lambda$ is a unitary character such that $\dim(B(\Lambda,\theta)^\I) = 1$, we use the function $B$ to generate a Hecke module $V_B$. We show that $V_B$ is irreducible and has a unique (up to a constant) vector which is Iwahori spherical, hence implying that $V_B$ is a $(\Lambda, \theta)$-Bessel model for the Steinberg representation. In case $\Lambda$ is not unitary (this can only happen when $L$ is split over $F$), we use the fact that the Steinberg representation, in the split case, is generic and actually show that \emph{any} generic, irreducible, admissible representation of $\GSp_4(F)$ has a split Bessel model. We believe that this result is known to the experts but since it is not available in the literature we present the proof in details. This is done is Sect.\ \ref{existence-section}, in particular Theorem \ref{existence-theorem}. Let us remark here that many of the results in Sect.\ \ref{bessel-model-section} involve a lot of computations. In the interest of brevity and to avoid distraction from the main point, we have sometimes illustrated the calculations for one case and left the other cases to the reader. For all the detailed computations, we refer the reader to a longer version of this manuscript \cite{P1} available on our homepage.

In the second half of the paper, we give the application of the explicit formula of the test vector. In Sect.\ \ref{non-arch-int-rep-section}, we obtain an integral representation of the $L$-function of the Steinberg representation twisted by any $\GL_2$ representation in Theorem \ref{steinberg-int-rep-L-fn}. In Theorems \ref{global-thm} and \ref{special values thm}, we obtain an integral representation for the global $L$-function and a special value result.

The author would like to thank Ralf Schmidt for all his help, in particular, for explaining how to obtain a Bessel model from a Whittaker model in the split case. The author would also like to thank Abhishek Saha for several fruitful discussions on this topic.

\section{Steinberg representation of $\GSp_4$}\label{pi-steinberg-setup}
\subsection*{Non-archimedean setup}
Let $F$ be a non-archimedean local field of characteristic zero.
Let $\OF$, $\p$, $\varpi$, $q$ be the ring of integers, prime ideal, uniformizer
and cardinality of the residue class field $\OF/\p$, respectively. Let us fix three
elements $a,b,c\in F$ such that $d := b^2-4ac \neq 0$. Let
\begin{equation}\label{Ldefeq}\renewcommand{\arraystretch}{1.2}
 L = \left\{
      \begin{array}{l@{\qquad\mbox{if }}l}
        F(\sqrt{d})&d \notin F^{\times2},\\
        F \oplus F&d \in F^{\times2}.
      \end{array}
    \right.
\end{equation}
In case $L = F \oplus F$, we consider $F$ diagonally embedded. If
$L$ is a field, we denote by $\bar x$ the Galois conjugate of
$x\in L$ over $F$. If $L=F\oplus F$, let $\overline{(x,y)}=(y,x)$.
In any case we let $N(x)=x\bar x$ and $\tr(x)=x+\bar x$.
We shall make the following {\bf assumptions}:
\begin{description}
 \item[(A1)] $a,b\in\OF$ and $c\in\OF^\times$.
 \item[(A2)] If $d \not\in F^{\times2}$, then $d$ is the generator of the discriminant of $L/F$.
  If $d \in F^{\times2}$, then $d \in \OF^{\times}$.
\end{description}

We set the Legendre symbol as follows,
\begin{equation}\label{legendresymboldefeq}
 \Big(\frac L{\p}\Big) := \left\{
                  \begin{array}{l@{\qquad\text{if }}l@{\qquad}l}
                    -1, &d \not\in F^{\times2},\:d \not\in \p&\mbox{(the inert case)}, \\
                    0, &d\not\in F^{\times2},\:d\in\p&\mbox{(the ramified case)}, \\
                    1, &d \in F^{\times2}&\mbox{(the split case)}.
                  \end{array}\right.
\end{equation}
If $L$ is a field, then let $\OF_L$ be its ring of integers.
If $L = F \oplus F$, then let $\OF_L = \OF \oplus \OF$.
Let $\varpi_L$ be the uniformizer of $\OF_L$ if $L$ is a field and set
$\varpi_L = (\varpi,1)$ if $L$ is not a field.
Note that, if $(\frac{L}{\p}) \neq -1$, then $N(\varpi_L) \in \varpi \OF^{\times}$. Let $\alpha \in \OF_L$ be defined by
\begin{equation}\label{alphadefeq}\renewcommand{\arraystretch}{2.1}
 \alpha:=\left\{\begin{array}{l@{\qquad\text{if }L}l}
 \displaystyle\frac{b+\sqrt{d}}{2c}&\text{ is a field},\\
 \displaystyle\Big(\frac{b+\sqrt{d}}{2c},\frac{b-\sqrt{d}}{2c}\Big)&=F\oplus F.
 \end{array}\right.
\end{equation}
We fix the following ideal in $\OF_L$,
\begin{equation}\label{ideal defn}\renewcommand{\arraystretch}{1.3}
 \P := \p\OF_L = \left\{
                  \begin{array}{l@{\qquad\text{if }}l}
                    \p_L & \big(\frac L{\p}\big) = -1,\\
                    \p_L^2 & \big(\frac L{\p}\big) = 0,\\
                    \p \oplus \p & \big(\frac L{\p}\big) = 1.
                  \end{array}
                \right.
\end{equation}
Here, $\p_L$ is the maximal ideal of $\OF_L$ when $L$ is a field
extension. Note that $\P$ is prime only if $\big(\frac
L\p\big)=-1$. We have $\P^n\cap\OF=\p^n$ for all $n\geq0$. Let us recall Lemma 3.1.1 of \cite{PS1}.
\begin{lemma}\label{quadraticextdisclemma-steinberg}
 Let notations be as above.
 \begin{enumerate}
  \item The
   elements $1$ and $\alpha$ constitute an integral basis of $L/F$.
  \item There exists no $x\in\OF$ such that $\alpha+x\in\P$.
\end{enumerate}
\end{lemma}

\subsection*{Steinberg representation}
Let us define the symplectic group $H=\GSp_4$ by
$$
 H(F) := \{g \in \GL_4(F) : \, ^tgJg = \mu_2(g)J,\:\mu_2(g)\in F^{\times} \},$$
where $J = \mat{}{1_2}{-1_2}{}$. The maximal compact subgroup is denoted by $K^H := \GSp_4(\OF)$. We define the Iwahori subgroup as follows,
\begin{equation}\label{Iwahori-gp-defn}
\I := \{g \in K^H:\:g\equiv\begin{bmatrix}\ast&0&\ast&\ast\\
 \ast&\ast&\ast&\ast\\0&0&\ast&\ast\\ 0&0&0&\ast\end{bmatrix}\pmod{\p} \}.
\end{equation}

Let $\Omega$ be an unramified, quadratic character of $F^\times$. Let $\pi$ be the Steinberg representation of $H(F)$, twisted by the character $\Omega$. This representation is denoted by $\Omega {\rm St}_{\GSp_4}$. Since we have assumed that $\Omega$ is quadratic, we see that $\pi$ has trivial central character. The Steinberg representation has the property that it is the only representation of $H(F)$ which has a unique (up to a constant) Iwahori fixed vector. The Iwahori Hecke algebra acts on the space of $\I$-invariant vectors. We will next describe the Iwahori Hecke algebra.

\subsection*{Iwahori Hecke algebra}
The Iwahori Hecke algebra ${\mathcal H}_\I$ of $H(F)$ is the convolution algebra of left and right $\I$-invariant functions on $H(F)$. We refer the reader to Sect.\ 2.1 of \cite{Sc} for details on the Iwahori Hecke algebra. Here, we state the two projection operators (projecting onto the Siegel and Klingen parabolic subgroups) and the Atkin Lehner involution. The unique (up to a constant) Iwahori fixed vector $v_0$ in $\pi$ is annihilated by the projection operators and is an eigenvector of the Atkin Lehner involution.
\begin{align}
& \sum\limits_{w \in \OF/\p} \pi(\begin{bmatrix}1&w\\&1\\&&1\\&&-w&1\end{bmatrix}) v_0 + \pi(s_1) v_0 = 0 \label{hecke-operator-steinberg-abstract-1} \\
& \pi(\eta_0) v_0 = \omega v_0                                                      \label{atkin-lehner-steinberg-abstract} \\
& \sum\limits_{y \in \OF/\p} \pi(\begin{bmatrix}1&&&\\&1&&\\y&&1&\\&&&1\end{bmatrix})  v_0 + \pi(s_2) v_0 = 0. \label{hecke-operator-steinberg-abstract-2}
\end{align}
Here
\begin{equation}\label{Weyl-group-elmt-defn}
s_1=\begin{bmatrix}&1\\1\\&&&1\\&&1\end{bmatrix},
 \qquad s_2=\begin{bmatrix}&&1\\&1\\-1\\&&&1\end{bmatrix}, \qquad \eta_0 = \begin{bmatrix}&&&1\\&&1&\\&\varpi&&\\ \varpi&&&\end{bmatrix} \mbox{ and } \omega = -\Omega(\varpi).
\end{equation}

\section{Existence and uniqueness of Bessel models for the Steinberg representation}\label{bessel-model-section}
Let us fix an additive character $\psi$ of $F$, with conductor $\OF$. Let $a,b \in \OF$ and $c \in \OF^\times$ be as in Sect.\ \ref{pi-steinberg-setup}, and set $S=\mat{a}{b/2}{b/2}{c}$. Then $\psi$ defines a character $\theta$ on $U(F)=\{\mat{1_2}{X}{}{1_2} : {}^tX = X\}$ by
\begin{equation}\label{thetadefeq-steinberg-pi}
 \theta(\mat{1}{X}{}{1})=\psi(\tr(SX)).
\end{equation}
Let
\begin{equation}\label{TRthetaeq1}
T(F) := \{ g \in \GL_2(F):\:^tgSg=\det(g)S\}.
\end{equation}
Set $\xi = \mat{b/2}{c}{-a}{b/2}$ and $F(\xi) = \{x+y\xi: x,y \in F\}$. Then, it can be checked that $T(F) = F(\xi)^\times$ and is isomorphic to $L^\times$, with the isomorphism given by
\begin{equation}\label{T-isomorphic-to-L}
\mat{x+\frac b2 y}{cy}{-ay}{x-\frac b2 y} \mapsto \left\{
                                                      \begin{array}{ll}
                                                        x + y \frac{\sqrt{d}}2, & \hbox{ if } L \mbox{ is a field;} \\
                                                        (x + y \frac{\sqrt{d}}2, x - y \frac{\sqrt{d}}2), & \hbox{ if } L = F \oplus F.
                                                      \end{array}
                                                    \right.
\end{equation}
We consider $T(F)$
as a subgroup of $H(F)$ via
$$
 T(F)\ni g\longmapsto\mat{g}{}{}{\det(g)\,^tg^{-1}}\in H(F).
$$
Let $R(F)=T(F)U(F)$. We call $R(F)$ the \emph{Bessel subgroup} of
$H(F)$ (with respect to the given data $a,b,c$).
Let $\Lambda$ be any character on $L^\times$ that is trivial on $F^\times$. We will consider $\Lambda$ as a character on $T(F)$. We have $\theta(t^{-1}ut)=\theta(u)$ for all $u\in U(F)$ and $t\in T(F)$. Hence the map $tu\mapsto\Lambda(t)\theta(u)$ defines a character of $R(F)$. We denote
this character by $\Lambda\otimes\theta$. 

As mentioned in the introduction, a linear functional $\beta : V \rightarrow\C$, satisfying $\beta(\pi(r)v) = (\Lambda \otimes \theta)(r) \beta(v)$ for any $r \in R(F), v \in V$, is called a $(\Lambda, \theta)$-Bessel functional for $\pi$. We say that $\pi$ has a $(\Lambda, \theta)$-Bessel model if $\pi$ is isomorphic to a subspace of smooth functions $B : H(F) \rightarrow \C$ satisfying
\begin{equation}\label{bessel-model-defining-condition}
B(tuh) = \Lambda(t) \theta(u) B(h) \qquad \mbox{ for all } t \in T(F), u \in U(F), h \in H(F).
\end{equation}
The existence of a non-zero $(\Lambda,\theta)$-Bessel functional for $\pi$ is equivalent to the existence of a non-trivial $(\Lambda,\theta)$-Bessel model for $\pi$. If $\pi$ has a non-zero $(\Lambda,\theta)$-Bessel functional $\beta$, then the space $\{B_v : v \in \pi, B_v(h):=\beta(\pi(h)v)\}$ gives a non-trivial $(\Lambda,\theta)$-Bessel model for $\pi$. Conversely, if $\pi$ has a non-trivial $(\Lambda,\theta)$-Bessel model $\{B_v : v \in \pi\}$ then the linear functional $\beta(v) := B_v(1)$ is a non-zero $(\Lambda,\theta)$-Bessel functional for $\pi$. We say that $v \in \pi$ is a test vector for a Bessel functional $\beta$ if $\beta(v) \neq 0$. Note that a vector $v \in \pi$ is a test vector for $\beta$ if and only if the corresponding function $B_v$ in the Bessel model satisfies $B_v(1) \neq 0$. 

Define the space $B(\Lambda, \theta)^\I$ of smooth functions $B$ on $H(F)$ which are right $\I$-invariant, satisfy (\ref{bessel-model-defining-condition}) and the following conditions, for any $h \in H(F)$, obtained from (\ref{hecke-operator-steinberg-abstract-1})-(\ref{hecke-operator-steinberg-abstract-2}),
\begin{align}
& \sum\limits_{w \in \OF/\p} B(h \begin{bmatrix}1&w\\&1\\&&1\\&&-w&1\end{bmatrix}) + B(h s_1) = 0, \label{hecke-operator-steinberg-1} \\
& B(h \eta_0) = \omega B(h),                                                         \label{atkin-lehner-steinberg} \\
& \sum\limits_{y \in \OF/\p} B(h \begin{bmatrix}1&&&\\&1&&\\y&&1&\\&&&1\end{bmatrix}) + B(h s_2) = 0. \label{hecke-operator-steinberg-2}
\end{align}
Our aim is to obtain the criteria for existence and uniqueness for $(\Lambda, \theta)$-Bessel models for $\pi$. Let us state the steps we take to obtain this. 
\begin{enumerate}
\item Since a function $B$ in $B(\Lambda, \theta)^\I$ is right $\I$-invariant and satisfies  (\ref{bessel-model-defining-condition}) we see that the values of $B$ are completely determined by its values on double coset representatives $R(F) \backslash H(F) / \I$. We obtain these representatives in Proposition \ref{disjoint-double-cosets-prop}.

\item In Proposition \ref{values-of-B-prop}, we use the $\I$-invariance of $B$ and (\ref{bessel-model-defining-condition})-(\ref{hecke-operator-steinberg-2}) to obtain necessary conditions  to be satisfied by the values of functions in $B(\Lambda, \theta)^\I$ on double coset representatives for $R(F) \backslash H(F) / \I$. This gives us ${\rm dim}(B(\Lambda, \theta)^\I) \leq 1$ in Corollary \ref{dim-B(lambda,theta)^I-leq-1-corollary}.

\item In Proposition \ref{well-defined-proposition}, we show that the function $B$ with the given values at double coset representatives for $R(F) \backslash H(F) / \I$ (obtained in Proposition \ref{values-of-B-prop}) is well-defined. We show that $B$  satisfies (\ref{hecke-operator-steinberg-1}), (\ref{atkin-lehner-steinberg}) and (\ref{hecke-operator-steinberg-2}) for all values of $h \in H(F)$ and obtain the criteria for ${\rm dim}(B(\Lambda, \theta)^\I) = 1$ in Theorem \ref{B(lambda-theta)^I-dimension-theorem}.

\item Suppose $\Lambda$ is such that ${\rm dim}(B(\Lambda, \theta)^\I) = 1$. If $\Lambda$ is unitary then we use $0 \neq B \in B(\Lambda, \theta)^\I$ to generate a Hecke module $V_B$. We define an inner product on $V_B$ and show in Proposition \ref{irreducibility-of-VB} that $V_B$ is irreducible and provides a $(\Lambda, \theta)$-Bessel model for $\pi$. If $\Lambda$ is not unitary (this can happen only if $L$ is a split extension of $F$), then we show that any irreducible, generic, admissible representation of $H(F)$ has a split $(\Lambda, \theta)$-Bessel model. Since $\pi$ is generic in the split case, we obtain in Theorem \ref{existence-theorem} the precise criteria for existence and uniqueness of a $(\Lambda, \theta)$-Bessel model for $\pi$.
\end{enumerate}

\subsection{Double coset decomposition}\label{double-coset-decomp-section}
From (3.4.2) of \cite{Fu}, we have the following disjoint double coset decomposition.
\begin{equation}\label{RFKHrepresentativeseq}
 H(F)=\bigsqcup_{l\in\Z}\bigsqcup_{m\geq0}R(F)h(l,m)K^H,
\end{equation}
where
\begin{equation}\label{hlmdefeq}
 h(l,m)=\begin{bmatrix}\varpi^{2m+l}\\&\varpi^{m+l}\\&&1\\&&&\varpi^m\end{bmatrix}.
\end{equation}

It follows from the Bruhat decomposition for $\SSp(4,\OF/\p)$ that
\begin{align}
 \label{newK2KHcosetseq41bA}K^H&=\I \sqcup\:\bigsqcup_{x\in\OF/\p}
  \begin{bmatrix}1\\x&1\\&&1&-x\\&&&1\end{bmatrix}s_1\I \sqcup\:\bigsqcup_{x\in\OF/\p}
  \begin{bmatrix}1&&x\\&1\\&&1&\\&&&1\end{bmatrix}s_2\I \sqcup \bigsqcup_{x,y\in\OF/\p}
  \begin{bmatrix}1\\x&1&&y\\&&1&-x\\&&&1\end{bmatrix}s_1s_2\I\\
 \label{newK2KHcosetseq44bA}&\sqcup\:\bigsqcup_{x,y\in\OF/\p}
  \begin{bmatrix}1&&x&y\\&1&y\\&&1\\&&&1\end{bmatrix}s_2s_1\I \sqcup\:\bigsqcup_{x,y,z\in\OF/\p}
  \begin{bmatrix}1&&&y\\x&1&y&xy+z\\&&1&-x\\&&&1\end{bmatrix}s_1s_2s_1\I \\
   \label{newK2KHcosetseq48bA} &\sqcup\:\bigsqcup_{x,y,z\in\OF/\p}
  \begin{bmatrix}1&&x&y\\&1&y&z\\&&1\\&&&1\end{bmatrix}s_2s_1s_2\I
\sqcup\:\bigsqcup_{w,x,y,z\in\OF/\p}
  \begin{bmatrix}1&&x&y\\w&1&wx+y&wy+z\\&&1&-w\\&&&1\end{bmatrix}s_1s_2s_1s_2\I.
\end{align}
Let $W = \{1, s_1, s_2, s_1 s_2, s_2 s_1, s_1 s_2 s_1, s_2 s_1 s_2, s_1 s_2 s_1 s_2\}$ be the Weyl group of $\SSp_4(F)$ and set
\begin{equation}\label{subsets-of-W}
W^{(1)} = \{1, s_2, s_2 s_1, s_2 s_1 s_2\}.
\end{equation}
Note that $W^{(1)}$ is the set of representatives for $\{1, s_1\}\backslash W$. Observing that $h(l,m) \begin{bmatrix}1&&\OF&\OF\\&1&\OF&\OF\\&&1&\\&&&1\end{bmatrix} h(l,m)^{-1}$ is contained in $R(F)$, we get a preliminary (non-disjoint) decomposition
\begin{equation}\label{prelim-decompA}
R(F) h(l,m) K^H = \bigcup_{s \in W^{(1)}, w \in \OF/\p} \Big(R(F) h(l,m) s \I \cup  R(F) h(l,m) W_w s_1 s \I \Big),
\end{equation}
where, for $w \in \OF/\p$, we set
$$W_w = \begin{bmatrix}1&&&\\w&1&&\\&&1&-w\\&&&1\end{bmatrix}.$$
The next lemma gives the condition under which the two double cosets of the form $R(F) h(l,m) s \I$ and  $R(F) h(l,m) W_w s_1 s \I$ are the same.
\begin{lemma}\label{coset-collapse-criterion-lemma}
For $w \in \OF/\p$ and $m \geq 0$, set $\beta_w^m := a \varpi^{2m}+b\varpi^mw+cw^2$. Let $s \in W^{(1)}$. Then $R(F) h(l,m) s \I = R(F) h(l,m) W_w s_1 s \I$ if and only if $\beta^m_w \in \OF^\times$.
\end{lemma}
{\bf Proof.} Suppose $\beta^m_w \in \OF^\times$. Take $y=\varpi^m, x=\varpi^mb/2+cw$ and set $g=\mat{x+\frac b2 y}{cy}{-ay}{x-\frac b2 y}$. Let $r= \mat{g}{}{}{\det(g)\,^tg^{-1}}$. Then
$$r h(l,m) = h(l,m) W_w s_1 k, \qquad \mbox{ where } k=\begin{bmatrix}-\beta_w^m&&&\\b\varpi^m+cw&c&&\\&&-c&b\varpi^m+cw\\&&&\beta_w^m\end{bmatrix} \in \I.$$
Note that for any $s \in W^{(1)}$, we have $s^{-1} k s \in \I$. Using $r h(l,m)s = h(l,m) W_w s_1 s (s^{-1} k
s)$, we obtain $R(F) h(l,m) s \I = R(F) h(l,m) W_w s_1 s \I$, as required.

Now we will prove the converse. For $s \in W^{(1)}$, let $r_s \in R(F)$ be such that $A_s= (h(l,m) W_w s_1s)^{-1} r_s h(l,m)s \in \I$. Then, it is possible to conclude that $\beta_w^m \in \OF^\times$. Let us illustrate the case $s = 1$ and $w \in (\OF/\p)^\times$. Let $r_1 = \mat{g}{}{}{\det(g)\,^tg^{-1}} \mat{1_2}{X}{}{1_2}$, where $X = \,^tX$ and $g = \mat{x+\frac b2 y}{cy}{-ay}{x-\frac b2 y}$.  Looking at the $(2,2)$ coefficient of $A_1$, we get $y = \varpi^m y'$, with $y' \in \OF^\times$. We see that the $(1,2)$ coefficient of  $A_1$ is equal to  $- \frac 1w ( (1,1)$ coefficient of $A_1 + \beta^m_w y')$ which lies in $\p$. Since diagonal elements of $A_1$ are in $\OF^\times$ and $w, y' \in \OF^\times$, we can conclude that $\beta_w^m \in \OF^\times$, as required. The other cases are done in a similar way (see \cite{P1}).
 \qed 

The next lemma describes for which $w \in \OF/\p$ we have $\beta^m_w \in \OF^\times$.
\begin{lemma}\label{beta-a-unit-lemma}
For $w \in \OF/\p$ and $m \geq 0$, set $\beta_w^m := a \varpi^{2m}+b\varpi^mw+cw^2$ as above.
\begin{enumerate}
\item If $m > 0$, then $\beta_w^m \in \OF^\times$ if and only if $w \in (\OF/\p)^\times$.

\item Let $m=0$.
\begin{enumerate}
\item If $\Big(\frac L{\p}\Big) = -1$, then $\beta_w^0 \in \OF^\times$ for every $w \in \OF/\p$.

\item Let $\Big(\frac L{\p}\Big) = 0$. Let $w_0$ be the unique element of $\OF/\p$ such that $\alpha + w_0 \in \p_L$, the prime ideal of $\OF_L$. Then $\beta_w^0 \in \OF^\times$ if and only if $w \neq w_0$. In case $\#(\OF/\p)$ is odd, one can take $w_0 = -b/(2c)$.

\item Let $\Big(\frac L{\p}\Big) = 1$. Then $\beta_w^0 \in \OF^\times$ if and only if $w \neq \frac{-b + \sqrt{d}}{2c}, \frac{-b - \sqrt{d}}{2c}$.
\end{enumerate}
\end{enumerate}
\end{lemma}
{\bf Proof.} Part $i)$ is clear. For the rest of the lemma, we need the following claim.

\underline{Claim}: We have $\beta_w^0 \in \OF^\times$ if and only if $\alpha + w \in \OF_L^\times$.

The claim follows from the identity
\begin{equation}\label{useful-identity-beta}
a+bw+cw^2 = -c(\alpha+w)(\bar\alpha + w) = -c N(\alpha+w).
\end{equation}
If $\Big(\frac L{\p}\Big) = -1$, then $\p_L = \P$ and Lemma \ref{quadraticextdisclemma-steinberg} ii) implies that $\alpha + w \in \OF_L^\times$ for all $w \in \OF/\p$. The claim gives ii)a) of the lemma. Let us now assume that $\Big(\frac L{\p}\Big) = 0$. In this case, the injective map $\iota : \OF \hookrightarrow \OF_L$ gives an isomorphism between the fields $\OF/\p \simeq \OF_L/\p_L$. Let $w_0 = -\iota^{-1}(\alpha)$ be the unique element in $\OF/\p$ such that $\alpha + w_0 \in \p_L$. In case $\#(\OF/\p)$ is odd, then one can take $w_0 = -b/(2c) \in \OF$ since $\sqrt{d} \in \p_L$. Then for any $w \in \OF/\p, w \neq w_0$, we have $\alpha + w \in \OF_L^\times$. Now, the claim gives ii)b) of the lemma. Next assume that $\Big(\frac L{\p}\Big) = 1$. Since $\sqrt{d} \in \OF^\times$ by assumption, we have $\alpha \not\in \P$. If $\alpha + w \not\in \OF_L^\times$ for some $w \in \OF$, then we have one of $(b\pm\sqrt{d})/(2c)+w$ lies in $\p$. Hence, we see that the only choices of $w=(w,w)$ such that $\alpha+w \not\in \OF_L^\times$ are $w = (-b \pm \sqrt{d})/(2c)$. Note that $\sqrt{d} \in \OF^\times$ implies that $(-b \pm \sqrt{d})/(2c)$ are not equal modulo $\p$. This completes the proof of the lemma. \qed

Note that, in the case $\Big(\frac L{\p}\Big) = 0$, (\ref{useful-identity-beta}) implies that $\beta_{w_0}^0 \in \p$ but $\beta_{w_0}^0 \not\in \p^2$ by Lemma \ref{quadraticextdisclemma-steinberg} ii).

Next, we will show the disjointness of all the relevant double cosets. Set
\begin{align*}
& A_{s,t} = (h(l,m) s)^{-1} r h(l,m) t \qquad \qquad \qquad \qquad \qquad  t \in W, s \in W - \{t, s_1 t\} \\
& A_{w,s,t} = (h(l,0) s)^{-1} r h(l,0) W_w s_1 t \qquad \qquad \qquad \quad \,\, w \in \OF/\p, s, t \in W^{(1)}, s \neq t \\
& A_{w,s,t}^* = (h(l,0) W_w s_1 s)^{-1} r h(l,0) W_w s_1 t \qquad \qquad \quad w \in \OF/\p, s, t \in W^{(1)}, s \neq t
\end{align*}
Notice that, for all of the matrices defined above, and any $r \in R(F)$, at least one of the diagonal entries is zero. This implies that none of the above matrices can be in $\I$ for any choice of $r \in R(F)$.
For $\Big(\frac L{\p}\Big) = 1$, set
$$A_{s,t}^* = (h(l,0) W_{\frac{-b+\sqrt{d}}{2c}} s_1 s)^{-1} r h(l,0) W_{\frac{-b-\sqrt{d}}{2c}} s_1 t \qquad \qquad \quad s, t \in W^{(1)}.$$
If $s \neq t$, then at least one of the diagonal entries of $A_{s,t}^*$ is zero, implying that it cannot be in $\I$ for any choice of $r \in R(F)$. If $s = t$ and $A_{s,t}^*$ is in $\I$, then we get
$$x - \frac{\sqrt{d}y}2 \in \OF^\times \mbox{ and } -\frac{\sqrt{d}}c (x - \frac{\sqrt{d}y}2)
 \in \p,$$
which is not possible, since $\sqrt{d} \in \OF^\times$. Hence, $A_{s,t}^*$ cannot be in $\I$ for any choice of $r \in R(F)$. We summarize in the following proposition.
\begin{proposition}\label{disjoint-double-cosets-prop}
Let $W$ be the Weyl group of $\SSp_4(F)$ and $W^{(1)} = \{1, s_2, s_2 s_1, s_2 s_1 s_2\}$. If $\Big( \frac L{\p}\Big) = 0$, let $w_0$ be the unique element of $\OF/\p$ such that $\alpha + w_0 \in \p_L$. If $\#(\OF/\p)$ is odd, then take $w_0 = -b/(2c)$. Then we have the following disjoint double coset decomposition.
\begin{equation}\label{disjoint-double-cosets-eqn}
R(F) h(l,m) K^H = \left\{
                    \begin{array}{ll}
                      \bigsqcup\limits_{s \in W} R(F) h(l,m) s \I , & \hbox{ if } m > 0; \\
                      \bigsqcup\limits_{s \in W^{(1)}} R(F) h(l,0) s \I, & \hbox{ if } m = 0, \Big( \frac L{\p}\Big) = -1; \\
                      \bigsqcup\limits_{s \in W^{(1)}} \Big(R(F) h(l,0) s \I \sqcup R(F) h(l,0) W_{w_0} s_1 s \I\Big), & \hbox{ if } m = 0, \Big( \frac L{\p}\Big) = 0; \\
                      \bigsqcup\limits_{s \in W^{(1)}} \Big(R(F) h(l,0) s \I \sqcup R(F) h(l,0) W_{\frac{-b + \sqrt{d}}{2c}} s_1 s & \\
                      \qquad \qquad \qquad  \sqcup R(F) h(l,0) W_{\frac{-b - \sqrt{d}}{2c}} s_1 s\I\Big) , & \hbox{ if } m = 0, \Big( \frac L{\p}\Big) = 1.
                    \end{array}
                  \right.
\end{equation}
\end{proposition}

\subsection{Necessary conditions for values of $B \in B(\Lambda, \theta)^\I$}\label{values-of-B-section}
In this section, we will obtain the necessary conditions on the values of $B \in B(\Lambda, \theta)^\I$ on the double coset representatives from Proposition \ref{disjoint-double-cosets-prop} using $\I$-invariance of $B$ and (\ref{bessel-model-defining-condition})-(\ref{hecke-operator-steinberg-2}).

\subsubsection*{Conductor of $\Lambda$:}
Let us define
\begin{equation}\label{c(Lambda)-defn}
c(\Lambda) = {\rm min}\{m \geq 0 : \Lambda |_{(1+\P^m)\cap \OF_L^\times} \equiv 1\}.
\end{equation}
Note that $(1+\P^m)\cap \OF_L^\times = 1+\P^m$ if $m \geq 1$ and $(1+\P^m)\cap \OF_L^\times = \OF_L^\times$ if $m = 0$. Also, $c(\Lambda)$ is the conductor of $\Lambda$ only if $\Big(\frac L{\p}\Big) = -1$. Let us set $c(\Lambda) = m_0$.  Since $\Lambda$ is trivial on $F^\times$, we see that $\Lambda|_{(\OF^\times + \P^{m_0}) \cap \OF_L^\times} \equiv 1$.

Let us make a few observations about $\Lambda$ and $c(\Lambda)$.
\begin{enumerate}
\item If $L$ is a field, then we have $L^\times = \langle \varpi_L \rangle . \OF_L^\times$. If $\Big(\frac L{\p}\Big) = -1$ and $m_0=0$, then we have that $\Lambda(\varpi_L) = 1$, since $\varpi_L \in \varpi \OF_L^\times$. In case $\Big(\frac L{\p}\Big) = 0$ and $m_0=0$, we see that $\Lambda(\varpi_L) = \pm 1$. In general, if $L$ is a field, we see that $\Lambda$ is a unitary character since $m_0$ is finite.
 \item If $L$ is not a field, then $L^\times = F^\times \oplus F^\times$ and $\Lambda((x,y)) = \Lambda_1(x) \Lambda_2(y)$, where $\Lambda_1, \Lambda_2$ are two characters of $F^\times$ satisfying $\Lambda_1.\Lambda_2 \equiv 1$. In this case, $m_0$ is the conductor of both $\Lambda_1, \Lambda_2$ and the character $\Lambda$ need not be unitary.
\end{enumerate}
In the next lemma, we will describe some coset representatives, which will be used in the evaluation of certain sums involving the character $\Lambda$.
\begin{lemma}\label{Lambda-conductor-coset-rep-lemma}
Let $m \geq 1$. A set of coset representatives for $((\OF^\times + \P^{m-1}) \cap \OF_L^\times)/(\OF^\times + \P^m)$ is given by
\begin{align*}
\{w + \alpha \varpi^{m-1} : w \in (\OF/\p)^\times \} \cup \{1\} & \quad \mbox{ if } m \geq 2 \\
\{w + \alpha : w \in \OF/\p \} \cup \{1\} & \quad \mbox{ if } m = 1, \Big(\frac{L}{\p}\Big) = -1 \\
\{w + \alpha : w \in \OF/\p, w \neq w_0 \} \cup \{1\} & \quad \mbox{ if } m = 1, \Big(\frac{L}{\p}\Big) = 0 \\
\{w + \alpha : w \in \OF/\p, w \neq (-b\pm\sqrt{d})/(2c) \} \cup \{1\} & \quad \mbox{ if } m = 1, \Big(\frac{L}{\p}\Big) = 1.
\end{align*}
In the case $ \Big(\frac{L}{\p}\Big) = 0$, the element $w_0$ is the unique element in $\OF/\p$ such that $w_0 + \alpha \notin \OF_L^\times$.
\end{lemma}
{\bf Proof.} Let $x+\alpha \varpi^{m-1}y \in (\OF^\times + \P^{m-1}) \cap \OF_L^\times$, with $x, y \in \OF$. If $m \geq 2$, then $x \in \OF^\times$. If $y \in \p$, then $x+\alpha \varpi^{m-1}y \in (\OF^\times + \P^m)$, and hence corresponds to the coset representative $1$. Now, let us assume that $y \in \OF^\times$. Then, using $y \in \OF^\times + \P^m$, we see that $x+\alpha \varpi^{m-1}y$ is equivalent to $x/y+\alpha \varpi^{m-1}$ modulo $(\OF^\times + \P^m)$. Note that $x/y+\alpha \varpi^{m-1} \in \OF_L^\times$ implies that, modulo $\p$, the element $x/y$ lies in
\begin{align}
 (\OF/\p)^\times \qquad \mbox{ if } m \geq 2, &
\quad \OF/\p \qquad \mbox{ if } m = 1, \Big(\frac{L}{\p}\Big) = -1 \nonumber\\
 \OF/\p - \{ w_0 \}  \quad \mbox{ if } m = 1, \Big(\frac{L}{\p}\Big) = 0, &
\quad \OF/\p - \{(-b\pm\sqrt{d})/(2c) \}  \quad \mbox{ if } m = 1, \Big(\frac{L}{\p}\Big) = 1. \label{possible-cosets}
\end{align}
This follows from the proof of Lemma \ref{beta-a-unit-lemma}. A calculation shows (see \cite{P1}) that if $w, w'$ are equivalent, modulo $\p$, to (not necessarily the same) elements in the sets defined in (\ref{possible-cosets}), then
$$w \equiv w' \pmod{\p} \quad \Leftrightarrow \quad (w+\alpha \varpi^{m-1})/(w'+\alpha \varpi^{m-1}) \in \OF^\times + \P^m.$$
This completes the proof of the lemma. \qed

Depending on the $c(\Lambda)$, certain values of $B$ have to be zero. This is obtained in the next lemma.

\begin{lemma}\label{conductor-of-lambda-B-vanishing-lemma}
\begin{enumerate}
\item Let $c(\Lambda) = m_0 \geq 2$. Then
\begin{equation}\label{B-vanishing-m-condition-equation}
B(h(l,m)) = 0 \qquad \mbox{ for all } m \leq m_0-2 \mbox{ and any } l.
\end{equation}

\item Let $c(\Lambda) = m_0 \geq 1$ and $\Big(\frac L{\p}\Big) = 1$. For $w = (-b\pm \sqrt{d})/(2c)$
$$B(h(l,0)W_w s_1) = 0 \mbox{ for all } l.$$

\item Let $c(\Lambda) = m_0 = 0, \Big(\frac L{\p}\Big) = 0$ and $\Lambda = \Omega \circ N_{L/F}$ . Then
$$B(h(l,0) W_{w_0} s_1 s_2) = 0 \mbox{ for all } l.$$

\item Let $c(\Lambda) = m_0 = 0$ and  $\Big(\frac L{\p}\Big) = -1$. Then
$$B(h(l,0)) = 0 \mbox{ for all } l.$$
\end{enumerate}
\end{lemma}
{\bf Proof.} Let us illustrate the proof of $i)$ here.
Let $m \leq m_0-2$. Let $1+x+\alpha y \in 1+\P^{m+1}, x,y \in \p^{m+1}$, such that $\Lambda(1 + x + \alpha y) \neq 1$. Let
$$k = \begin{bmatrix}c(1+x)+by&cy\varpi^{-m}&&\\-ay\varpi^m&c(1+x)&&\\&&c(1+x)&ay\varpi^m\\&&-cy\varpi^{-m}&c(1+x)+by\end{bmatrix} \in \I.$$
Then
\begin{align*}
B(h(l,m)) &= B(h(l,m) k) = B(\begin{bmatrix}c(1+x)+by&cy&&\\-ay&c(1+x)&&\\&&c(1+x)&ay\\&&-cy&c(1+x)+by\end{bmatrix} h(l,m)) \\
&= \Lambda(1+x+\alpha y) B(h(l,m)),
\end{align*}
which implies that $B(h(l,m)) = 0$, as required. The other cases are computed in a similar manner (see \cite{P1}). \qed

From Lemmas \ref{Lambda-conductor-coset-rep-lemma} and \ref{conductor-of-lambda-B-vanishing-lemma}(i), we obtain the following information on certain character sums involving $\Lambda$.
\begin{lemma}\label{Lambda-character-sum-lemma} For any $l$, we have
\begin{align}
\sum\limits_{w \in (\OF/\p)^\times} \Lambda(w + \alpha \varpi^m) B(h(l,m)) + B(h(l,m)) &= \left\{
                                                                                             \begin{array}{ll}
                                                                                               0, & \hbox{ if } m < m_0; \\
                                                                                               q B(h(l,m)), & \hbox{ if } m \geq m_0.
                                                                                             \end{array}
                                                                                           \right. \mbox{ if } m > 0. \label{lambda-m>0-sum}\\
\sum\limits_{w \in \OF/\p} \Lambda(w + \alpha) B(h(l,0)) + B(h(l,0)) &= \left\{
                                                                                             \begin{array}{ll}
                                                                                               0, & \hbox{ if } m_0 \geq 1; \\
                                                                                               (q+1) B(h(l,0)), & \hbox{ if } m_0=0.
                                                                                             \end{array}
                                                                                           \right. \mbox{ if } \Big(\frac{L}{\p}\Big) = -1 \label{lambda-m=0-inert-sum}\\
 \sum\limits_{\substack{w \in \OF/\p\\ w \neq w_0}} \Lambda(w + \alpha) B(h(l,0)) + B(h(l,0)) &= \left\{
                                                                                             \begin{array}{ll}
                                                                                               0, & \hbox{ if } m_0 \geq 1; \\
                                                                                               q B(h(l,0)), & \hbox{ if } m_0=0.
                                                                                             \end{array}
                                                                                           \right. \mbox{ if } \Big(\frac{L}{\p}\Big) = 0 \label{lambda-m=0-ramified-sum}\\
\sum\limits_{\substack{w \in \OF/\p\\ w \neq \frac{-b\pm\sqrt{d}}{2c}}} \Lambda(w + \alpha) B(h(l,0)) + B(h(l,0)) &= \left\{
                                                                                             \begin{array}{ll}
                                                                                               0, & \hbox{ if } m_0 \geq 1; \\
                                                                                               (q-1) B(h(l,0)), & \hbox{ if } m_0=0.
                                                                                             \end{array}
                                                                                           \right. \mbox{ if } \Big(\frac{L}{\p}\Big) = 1\label{lambda-m=0-split-sum}
\end{align}
\end{lemma}

\subsubsection*{Conductor of $\psi$}
Since the conductor of $\psi$ is $\OF$, we obtain the following further vanishing conditions on the values of $B$.
\begin{lemma}\label{B-vanishing-conditions}
\begin{enumerate}
\item If $s \in \{1, s_1, s_2, s_2 s_1\}$ and $m \geq 0$, then
$$B(h(l,m)s) = 0 \qquad \mbox{ if } l < 0.$$

\item If $s \in \{s_1 s_2, s_1 s_2 s_1, s_2 s_1 s_2, s_1 s_2 s_1 s_2 \}$ and $m \geq 0$, then
$$B(h(l,m) s) = 0 \qquad \mbox{ if } l < -1.$$

\item If $w \in \OF$, then
$$B(h(l,0) W_w s_1) = 0 \qquad \mbox{ if } l < 0.$$

\item If $w \in \OF$ and $s \in \{s_1 s_2, s_1 s_2 s_1, s_1 s_2 s_1 s_2 \}$, then
$$B(h(l,0) W_w s) = 0 \qquad \mbox{ if } l < -1.$$

\item If $\Big(\frac L{\p}\Big) = 1$ and $w = (-b\pm \sqrt{d})/(2c)$, then
$$B(h(-1,0)W_w s_1 s_2) = 0.$$
\end{enumerate}
\end{lemma}
{\bf Proof.} Let us illustrate the proof of $i)$. For any $\epsilon \in \OF^\times$, set
$$k_s^\epsilon = \begin{bmatrix}1\\&1&&\epsilon\\&&1\\&&&1\end{bmatrix} \mbox{ if } s = 1, s_2 \mbox{ and } k_s^\epsilon = \begin{bmatrix}1&&\epsilon\\&1\\&&1\\&&&1\end{bmatrix} \mbox{ if } s = s_1, s_2 s_1.$$
Then, for $s \in \{1, s_1, s_2, s_2 s_1\}$ and $\epsilon \in \OF^\times$, we obtain
$$
B(h(l,m)s) = B(h(l,m) s k_s^\epsilon) = B(\begin{bmatrix}1&&&\\&1&&\epsilon \varpi^l\\&&1&\\&&&1\end{bmatrix} h(l,m)s)
= \psi(c\epsilon \varpi^l) B(h(l,m)s).
$$
Since the conductor of $\psi$ is $\OF$, we conclude that $B(h(l,m)s) = 0$ if $l < 0$. The other cases are computed in a similar manner (see \cite{P1}). \qed

\subsubsection*{Values of $B$ using (\ref{hecke-operator-steinberg-1})}
Substituting $h = h(l,m) s_1$ in (\ref{hecke-operator-steinberg-1}) and using Lemmas \ref{coset-collapse-criterion-lemma}, \ref{beta-a-unit-lemma} and \ref{Lambda-character-sum-lemma}, we get for any $l$
\begin{align}
B(h(l,m) s_1) &= \left\{
                                                                                             \begin{array}{ll}
                                                                                               0, & \hbox{ if } m < m_0; \\
                                                                                               -q B(h(l,m)), & \hbox{ if } m \geq m_0,
                                                                                             \end{array}
                                                                                           \right. \qquad \mbox{ if } m > 0. \label{B(h(l,m)s1)-m>0-value}\\
B(h(l,0) W_{w_0} s_1) &= \left\{
                                                                                             \begin{array}{ll}
                                                                                               0, & \hbox{ if } m_0 \geq 1; \\
                                                                                               -q B(h(l,0)), & \hbox{ if } m_0=0.
                                                                                             \end{array}
                                                                                           \right. \label{B(h(l,m)s1)-m=0-ramified-value}\\
B(h(l,0) W_{\frac{-b+\sqrt{d}}{2c}} s_1) + B(h(l,0) W_{\frac{-b-\sqrt{d}}{2c}} s_1) &= -(q-1) B(h(l,0)) \qquad \mbox{ if } m_0=0. \label{B(h(l,m)s1)-m=0-split-value}
\end{align}

Substituting $h = h(l,m) s_2 s_1$ in (\ref{hecke-operator-steinberg-1}) and using that the conductor of $\psi$ is $\OF$, we get for any $l, m$
\begin{equation}\label{B(h(l,m)s2s1)-value}
B(h(l,m) s_2 s_1) = -\frac 1{q} B(h(l,m)s_2).
\end{equation}

Substituting $h = h(l,m) s_1 s_2 s_1$ in (\ref{hecke-operator-steinberg-1}) and using that the conductor of $\psi$ is $\OF$, we get for any $m > 0$ and $l$
\begin{equation}\label{B(h(l,m)s1s2s1)-m>0-in-terms-of-s1s2-value}
B(h(l,m) s_1 s_2 s_1) = -\frac 1{q} B(h(l,m) s_1 s_2).
\end{equation}

Let $\Big(\frac L{\p}\Big) = 0$. Substituting $h = h(-1,0) W_{w_0} s_1 s_2 s_1$ in (\ref{hecke-operator-steinberg-1}) and using that the conductor of $\psi$ is $\OF$ and $b + 2c w_0 \in \p$, we get
\begin{equation}
B(h(-1,0) W_{w_0} s_1 s_2 s_1) = -\frac 1q B(h(-1,0) W_{w_0} s_1 s_2).
\end{equation}

Let $\Big(\frac L{\p}\Big) = 1$ and $w = (-b \pm \sqrt{d})/(2c)$. Substituting $h = h(l,0) W_w s_1 s_2 s_1$ in (\ref{hecke-operator-steinberg-1}) and using that the conductor of $\psi$ is $\OF$ and $\sqrt{d} \in \OF^\times$, we get for $l \neq -1$
\begin{equation}\label{B-value-s1s2s1-in-terms-of-s1s2}
B(h(l,m) W_w s_1 s_2 s_1) = - \frac 1q B(h(l,m) W_w s_1 s_2).
\end{equation}

\subsubsection*{Values of $B$ using (\ref{hecke-operator-steinberg-2})}
Substituting $h = h(l,m) s_2$ in (\ref{hecke-operator-steinberg-2}) and using that the conductor of $\psi$ is $\OF$, we get for any $l, m$
\begin{equation}\label{B(h(l,m)s2)-value}
B(h(l,m) s_2) = -\frac 1{q} B(h(l,m)).
\end{equation}

Substituting $h = h(l,m) s_2 s_1 s_2$ in (\ref{hecke-operator-steinberg-2}) and using that the conductor of $\psi$ is $\OF$, we get for $l \neq -1$
\begin{equation}\label{B-value-s2s1s2-in-terms-of-s2s1}
B(h(l,m) s_2 s_1 s_2) = - \frac 1q B(h(l,m) s_2 s_1).
\end{equation}

Let $w = 0$ if $m > 0$, $w = w_0$ if $m = 0, \Big(\frac L{\p}\Big) = 0$ and $w = (-b \pm \sqrt{d})/(2c)$ if $m = 0, \Big(\frac L{\p}\Big) = 1$. Substituting $h = h(l,m) W_w s_1 s_2$ in (\ref{hecke-operator-steinberg-2}) and using that the conductor of $\psi$ is $\OF$, we get for $l \neq -1$
\begin{equation}\label{B-value-s1s2-in-terms-of-s1}
B(h(l,m) W_w s_1 s_2) = - \frac 1q B(h(l,m) W_w s_1).
\end{equation}
Substituting $h = h(l,m) W_w s_1 s_2 s_1 s_2$ in (\ref{hecke-operator-steinberg-2}) and using that the conductor of $\psi$ is $\OF$, we get for all $l,m$
\begin{equation}\label{B-value-s1s2s1s2-lm-in-terms-of-s1s2s1}
B(h(l,m) W_w s_1 s_2 s_1 s_2) = - \frac 1q B(h(l,m) W_w s_1 s_2 s_1).
\end{equation}

\subsubsection*{Values of $B$ using (\ref{atkin-lehner-steinberg})}
    For any $l,m,w$ we have the matrix identities
\begin{align}
h(l,m) s_2 s_1 \eta_0 &= h(l-1,m+1) s_1 s_2 s_1 \begin{bmatrix}1&&&\\&-1&&\\&&-1&\\&&&1\end{bmatrix} \\
h(l,m) W_w s_1 s_2 s_1 s_2 \eta_0 &= h(l+1,m) W_w s_1 \begin{bmatrix}1&&&\\&1&&\\&&-1&\\&&&-1\end{bmatrix} \\
h(l,m) s_2 s_1 s_2 \eta_0 &= h(l+1,m) \begin{bmatrix}1&&&\\&1&&\\&&-1&\\&&&-1\end{bmatrix}. \label{eta-3rd-eqn}
\end{align}
Hence, by (\ref{atkin-lehner-steinberg}), we have
\begin{align}
B(h(l,m) s_2 s_1) &= \omega B(h(l-1,m+1) s_1 s_2 s_1), \label{B-value-lm-s2s1-in-terms-of-l-1-m+1-s1s2s1}\\
B(h(l,m)W_w s_1 s_2 s_1 s_2) &= \omega B(h(l+1,m) W_w s_1), \label{B-value-s2s1s2s1-lm-in-terms-of-l+1m-s1}\\
B(h(l,m)s_2 s_1 s_2) &= \omega B(h(l+1,m)). \label{B-value-s2s1s2-lm-in-terms-of-l+1m}
\end{align}
Using (\ref{eta-3rd-eqn}) we see that
$$B(h(l,0)W_{\frac{-b+\sqrt{d}}{2c}} s_1 s_2) = \omega B(h(l,0)W_{\frac{-b+\sqrt{d}}{2c}} s_1 s_2 \eta_0) = \omega B(h(l,0)W_{\frac{-b+\sqrt{d}}{2c}} \begin{bmatrix}1\\&\varpi\\&&\varpi\\&&&1\end{bmatrix} s_2).$$
Let $x = \sqrt{d}/2 + \varpi, y = 1, g = \mat{x+by/2}{cy}{-ay}{x-by/2}$ and set $r = \mat{g}{}{}{\det(g)\,^tg^{-1}}$. Then we have the matrix identity
$$r h(l,0) W_{\frac{-b-\sqrt{d}}{2c}} s_1 s_2 = h(l,0)W_{\frac{-b+\sqrt{d}}{2c}} \begin{bmatrix}1\\&\varpi\\&&\varpi\\&&&1\end{bmatrix} s_2 k, \mbox{ with } k = \begin{bmatrix}\frac{\sqrt{d}}c&&&-1\\&-\frac{\sqrt{d}}c&1\\&\varpi&c\\-\varpi&&&-c\end{bmatrix} \in \I.$$
This gives us
\begin{equation}\label{B-value-split-w1s1s2-in-terms-of-w2s1s2}
B(h(l,0)W_{\frac{-b+\sqrt{d}}{2c}} s_1 s_2) = \omega \Lambda((\sqrt{d}+\varpi, \varpi)) B(h(l,0) W_{\frac{-b-\sqrt{d}}{2c}} s_1 s_2).
\end{equation}

\subsubsection*{Summary}
Using  (\ref{B(h(l,m)s2s1)-value}), (\ref{B(h(l,m)s2)-value}), (\ref{B-value-s2s1s2-in-terms-of-s2s1}) and (\ref{B-value-s2s1s2-lm-in-terms-of-l+1m}) we get for $l, m \geq 0$
\begin{equation}\label{B-value-l+1-in-terms-of-l}
B(h(l+1,m)) = -\frac{\omega}{q^3} B(h(l,m)).
\end{equation}
Using (\ref{B(h(l,m)s1)-m>0-value}), (\ref{B(h(l,m)s2s1)-value}), (\ref{B(h(l,m)s1s2s1)-m>0-in-terms-of-s1s2-value}), (\ref{B(h(l,m)s2)-value}), (\ref{B-value-s1s2-in-terms-of-s1}), (\ref{B-value-lm-s2s1-in-terms-of-l-1-m+1-s1s2s1}) and (\ref{B-value-l+1-in-terms-of-l}), we get for $l \geq 0, m \geq m_0 - 1$
\begin{equation}\label{B-value-m+1-in-terms-of-m}
B(h(l,m+1)) = \frac 1{q^4} B(h(l,m)).
\end{equation}
Hence, we conclude that
\begin{equation}\label{B(h(l,m))-values}
B(h(l,m)) = \left\{
              \begin{array}{ll}
                0, & \hbox{ if } l \leq -1 \mbox{ or } 0 \leq m \leq m_0-2;\\
                q^{-4(m-m_0+1)}(-\omega q^{-3})^l B(h(0, m_0-1)), & \hbox{ if } l \geq 0 \mbox{ and } m \geq m_0 - 1 > 0; \\
                q^{-4m}(-\omega q^{-3})^l B(1), & \hbox{ if } l \geq 0 \mbox{ and } m \geq m_0 = 0, 1.
              \end{array}
            \right.
\end{equation}
Let $\Big(\frac{L}{\p}\Big) = 1$ and $w = (-b \pm \sqrt{d})/(2c)$. Using (\ref{B-value-s1s2s1-in-terms-of-s1s2}), (\ref{B-value-s1s2-in-terms-of-s1}), (\ref{B-value-s1s2s1s2-lm-in-terms-of-s1s2s1}) and (\ref{B-value-s2s1s2s1-lm-in-terms-of-l+1m-s1}), we get for $l \geq 0$
\begin{equation}
B(h(l+1,0) W_w s_1) = -\frac{\omega}{q^3} B(h(l,0) W_w s_1),
\end{equation}
which gives us
\begin{equation}
B(h(l,0) W_w s_1) = (-\omega q^{-3})^l B(W_w s_1).
\end{equation}
In addition, if $m_0 = 0$ and $\omega \Lambda((1,\varpi)) = -1$, using (\ref{B(h(l,m)s1)-m=0-split-value}), (\ref{B-value-s1s2-in-terms-of-s1}) and (\ref{B-value-split-w1s1s2-in-terms-of-w2s1s2}), we get for all $l \geq 0$
\begin{equation}
B(h(l,0)) = 0.
\end{equation}

Summarizing the calculations of the values of $B$, we obtain
\begin{proposition}\label{values-of-B-prop}
Let $c(\Lambda) = m_0$. For $l, m \in \Z, m \geq 0$, let us set
$$A_{l,m} := \left\{\begin{array}{ll}
                          q^{-4(m-m_0+1)}(-\omega q^{-3})^l, & \mbox{ if } m_0 \geq 1;\\
                          q^{-4m}(-\omega q^{-3})^l, & \mbox{ if } m_0 = 0.
                          \end{array}
                          \right. \qquad 
                  C_{m_0} := \left\{\begin{array}{ll}
                          B(h(0, m_0-1)), & \mbox{ if } m_0 \geq 1;\\
                          B(1), & \mbox{ if } m_0 = 0.
                          \end{array}
                          \right. $$
We have the following necessary conditions on the values of $B \in B(\Lambda, \theta)^\I$.
\begin{enumerate}
\item For $m \geq 0$ and any $m_0$,
\begin{enumerate}
\item $$B(h(l,m)) = \left\{
                      \begin{array}{ll}
                        0, & \hbox{ if } l \leq -1 \mbox{ or } m \leq m_0-2; \\
                        A_{l,m} C_{m_0}, & \hbox{ if } l \geq 0 \mbox{ and } m \geq m_0-1.
                      \end{array}
                    \right. $$

\item $$B(h(l,m)s_2) = \left\{
                      \begin{array}{ll}
                        0, & \hbox{ if } l \leq -1 \mbox{ or } m \leq m_0-2; \\
                        -\frac 1qA_{l,m} C_{m_0}, & \hbox{ if } l \geq 0 \mbox{ and } m \geq m_0-1.
                      \end{array}
                    \right. $$
\item $$B(h(l,m)s_2s_1) = \left\{
                      \begin{array}{ll}
                        0, & \hbox{ if } l \leq -1 \mbox{ or } m \leq m_0-2; \\
                        \frac 1{q^2}A_{l,m} C_{m_0}, & \hbox{ if } l \geq 0 \mbox{ and } m \geq m_0-1.
                      \end{array}
                    \right. $$

\item $$B(h(l,m)s_2s_1s_2) = \left\{
                            \begin{array}{ll}
                              0, & \hbox{ if } l \leq -2 \mbox{ or } m \leq m_0-2; \\
                              \omega A_{0,m} C_{m_0}, & \hbox{ if } l = -1 \mbox{ and } m \geq m_0-1; \\
                              -\frac 1{q^3}A_{l,m} C_{m_0}, & \hbox{ if } l \geq 0 \mbox{ and } m \geq m_0-1.
                            \end{array}
                          \right.$$
\end{enumerate}
\item For $m > 0$ and any $m_0$,
\begin{enumerate}
\item $$B(h(l,m)s_1) = \left\{
                      \begin{array}{ll}
                        0, & \hbox{ if } l \leq -1 \mbox{ or } m \leq m_0-1; \\
                        -qA_{l,m} C_{m_0}, & \hbox{ if } l \geq 0 \mbox{ and } m \geq m_0.
                      \end{array}
                    \right. $$

\item $$B(h(l,m)s_1s_2) = \left\{
                            \begin{array}{ll}
                              0, & \hbox{ if } l \leq -2 \mbox{ or } m \leq m_0-1; \\
                              -\omega q^3 A_{0,m} C_{m_0}, & \hbox{ if } l = -1 \mbox{ and } m \geq m_0; \\
                              A_{l,m} C_{m_0}, & \hbox{ if } l \geq 0 \mbox{ and } m \geq m_0.
                            \end{array}
                          \right.$$

\item $$B(h(l,m)s_1s_2s_1) = \left\{
                            \begin{array}{ll}
                              0, & \hbox{ if } l \leq -2 \mbox{ or } m \leq m_0-1; \\
                              \omega q^2 A_{0,m} C_{m_0}, & \hbox{ if } l = -1 \mbox{ and } m \geq m_0; \\
                              -\frac 1qA_{l,m} C_{m_0}, & \hbox{ if } l \geq 0 \mbox{ and } m \geq m_0.
                            \end{array}
                          \right.$$

\item $$B(h(l,m)s_1s_2s_1s_2) = \left\{
                            \begin{array}{ll}
                              0, & \hbox{ if } l \leq -2 \mbox{ or } m \leq m_0-1; \\
                              -\omega q A_{0,m} C_{m_0}, & \hbox{ if } l = -1 \mbox{ and } m \geq m_0; \\
                              \frac 1{q^2}A_{l,m} C_{m_0}, & \hbox{ if } l \geq 0 \mbox{ and } m \geq m_0.
                            \end{array}
                          \right.$$
\end{enumerate}

\item Let $m_0 \geq 1$.
\begin{enumerate}
\item If $\Big(\frac L{\p}\Big) = 0$ and $s \in \{1, s_2, s_2s_1, s_2 s_1 s_2\}$, then, for all $l$,
$$B(h(l,0)W_{w_0} s_1 s) = 0.$$

\item If $\Big(\frac L{\p}\Big) = 1, s \in \{1, s_2, s_2s_1, s_2 s_1 s_2\}$ and $w = \frac{-b \pm \sqrt{d}}{2c}$, then, for all $l$,
$$B(h(l,0)W_w s_1 s) = 0.$$
\end{enumerate}

\item Let $m_0 = 0$.
\begin{enumerate}
\item If $\Big(\frac L{\p}\Big) = -1$ then
$$C_0 = 0.$$

\item Suppose $\Big(\frac L{\p}\Big) = 0$, then
\begin{enumerate}
\item $$B(h(l,0)W_{w_0}s_1) = \left\{
                                \begin{array}{ll}
                                  0, & \hbox{ if } l \leq -1; \\
                                  -q A_{l,0} C_0, & \hbox{ if } l \geq 0.
                                \end{array}
                              \right.$$

\item $$B(h(l,0)W_{w_0}s_1 s_2) = \left\{
                                \begin{array}{ll}
                                  0, & \hbox{ if } l \leq -2; \\
-\omega q^3 C_0, & \hbox{ if } l = -1; \\
                                  A_{l,0} C_0, & \hbox{ if } l \geq 0.
                                \end{array}
                              \right.$$

\item $$B(h(l,0)W_{w_0}s_1 s_2 s_1) = \left\{
                                \begin{array}{ll}
                                  0, & \hbox{ if } l \leq -2; \\
                                  \omega q^2 A_{l+1,0} C_0, & \hbox{ if } l \geq -1.
                                \end{array}
                              \right.$$

\item $$B(h(l,0)W_{w_0}s_1 s_2 s_1 s_2) = \left\{
                                \begin{array}{ll}
                                  0, & \hbox{ if } l \leq -2; \\
                                  -\omega q A_{l+1,0} C_0, & \hbox{ if } l \geq -1.
                                \end{array}
                              \right.$$

\end{enumerate}

\item Suppose $\Big(\frac L{\p}\Big) = 0$ and $\Lambda = \Omega \circ N_{L/F}$, then
$$C_0 = 0.$$

\item Suppose $\Big(\frac L{\p}\Big) = 1$. Then for $s\in \{1, s_2, s_2 s_1, s_2 s_1 s_2\}$
 $$B(h(l,0)W_{\frac{-b - \sqrt{d}}{2c}} s_1 s)  = \frac{1}{\omega \Lambda((1,\varpi))} B(h(l,0)W_{\frac{-b + \sqrt{d}}{2c}} s_1 s).$$

\item Suppose $\Big(\frac L{\p}\Big) = 1$ and $\omega \Lambda((1,\varpi)) = -1$.
\begin{enumerate}
\item $$C_0 = 0.$$

\item $$B(h(l,0)W_{\frac{-b + \sqrt{d}}{2c}} s_1) = \left\{
                                \begin{array}{ll}
                                  0, & \hbox{ if } l \leq -1; \\
                                  A_{l,0} B(W_{\frac{-b + \sqrt{d}}{2c}} s_1), & \hbox{ if } l \geq 0.
                                \end{array}
                              \right.$$

\item $$B(h(l,0)W_{\frac{-b + \sqrt{d}}{2c}} s_1 s_2) = \left\{
                                \begin{array}{ll}
                                  0, & \hbox{ if } l \leq -1; \\
                                  -\frac 1qA_{l,0} B(W_{\frac{-b + \sqrt{d}}{2c}} s_1), & \hbox{ if } l \geq 0.
                                \end{array}
                              \right.$$

\item $$B(h(l,0)W_{\frac{-b + \sqrt{d}}{2c}} s_1 s_2 s_1) = \left\{
                                \begin{array}{ll}
                                  0, & \hbox{ if } l \leq -2; \\
                                  - \omega q A_{l+1,0} B(W_{\frac{-b + \sqrt{d}}{2c}} s_1), & \hbox{ if } l \geq -1.
                                \end{array}
                              \right.$$

\item $$B(h(l,0)W_{\frac{-b + \sqrt{d}}{2c}} s_1 s_2 s_1 s_2) = \left\{
                                \begin{array}{ll}
                                  0, & \hbox{ if } l \leq -2; \\
                                  \omega A_{l+1,0} B(W_{\frac{-b + \sqrt{d}}{2c}} s_1), & \hbox{ if } l \geq -1.
                                \end{array}
                              \right.$$

\end{enumerate}
\item Suppose $\Big(\frac L{\p}\Big) = 1$ and $\omega \Lambda((1,\varpi)) \neq -1$.

\begin{enumerate}
\item $$B(h(l,0) W_{\frac{-b + \sqrt{d}}{2c}} s_1) = \left\{
                                                                         \begin{array}{ll}
                                                                           0, & \hbox{ if } l \leq -1;\\
                                                                           -\frac{q-1}{1+\omega \Lambda((1,\varpi))} A_{l,0} C_0, & \hbox{ if } l \geq 0.
                                                                         \end{array}
                                                                       \right.$$

\item $$B(h(l,0) W_{\frac{-b + \sqrt{d}}{2c}} s_1 s_2) = \left\{
                                                                         \begin{array}{ll}
                                                                           0, & \hbox{ if } l \leq -1;\\
                                                                           \frac{q-1}{q(1+\omega \Lambda((1,\varpi)))} A_{l,0} C_0, & \hbox{ if } l \geq 0.
                                                                         \end{array}
                                                                       \right.$$

\item $$B(h(l,0) W_{\frac{-b + \sqrt{d}}{2c}} s_1 s_2 s_1) = \left\{
                                                                         \begin{array}{ll}
                                                                           0, & \hbox{ if } l \leq -2;\\
                                                                           \frac{\omega q(q-1)}{1+\omega \Lambda((1,\varpi))} A_{l+1,0} C_0, & \hbox{ if } l \geq -1.
                                                                         \end{array}
                                                                       \right.$$

\item $$B(h(l,0) W_{\frac{-b + \sqrt{d}}{2c}} s_1 s_2 s_1 s_2) = \left\{
                                                                         \begin{array}{ll}
                                                                           0, & \hbox{ if } l \leq -2;\\
                                                                           -\frac{\omega (q-1)}{1+\omega \Lambda((1,\varpi))} A_{l+1,0} C_0, & \hbox{ if } l \geq -1.
                                                                         \end{array}
                                                                       \right.$$
\end{enumerate}
\end{enumerate}

\end{enumerate}

\end{proposition}

The above proposition immediately gives us the following corollary.
\begin{corollary}\label{dim-B(lambda,theta)^I-leq-1-corollary}
For any character $\Lambda$, we have
\begin{equation}\label{dim-B(lambda,theta)^I-leq-1-eqn}
{\rm dim}\big(B(\Lambda, \theta)^\I\big) \leq 1.
\end{equation}
\end{corollary}

\subsection{Well-definedness of $B$}
In this section, we will show that a function $B$ on $H(F)$, which is right $\I$-invariant, satisfies (\ref{bessel-model-defining-condition}) and with values on the double coset representatives of $R(F) \backslash H(F) / \I$ given by Proposition \ref{values-of-B-prop}, is well defined. Hence, we have to show that
$$r_1 s k_1 = r_2 s k_2 \Rightarrow B(r_1 s k_1) = B(r_2 s k_2)$$
for $r_1, r_2 \in R(F), k_1, k_2 \in \I$ and any double coset representative $s$. This is
obtained in the following proposition.
\begin{proposition}\label{well-defined-proposition}
Let $s$ be any double coset representative from Proposition \ref{disjoint-double-cosets-prop} and the values $B(s)$ be as in Proposition \ref{values-of-B-prop}. Let $t \in T(F), u \in U(F)$ such that $s^{-1} t u s \in \I$. Then
$$\Lambda(t) \theta(u) = 1 \mbox{ or } B(s) = 0.$$
\end{proposition}
{\bf Proof.} Let $t = \mat{g}{}{}{\det(g)\,^tg^{-1}}$ and $u = \mat{1}{X}{}{1}$, with $g = \mat{x+by/2}{cy}{-ay}{x-by/2}, X = {}^tX$. First let  $s = h(l,m) s'$ with $s' \in W$. If $s' \in \{s_1, s_1 s_2, s_1 s_2 s_1, s_1 s_2 s_1 s_2\}$ then we only consider $m > 0$. Observe that $x+y\frac{\sqrt{d}}2 = x - \frac {by}2 + cy\alpha$. (In the split case, we consider the same identity with $(x+y\frac{\sqrt{d}}2, x-y\frac{\sqrt{d}}2)$). Let us assume $s^{-1} t u s \in \I$. For any $s'$, we see that $x \pm by/2 \in \OF^\times$. If $s' \in \{1, s_2, s_2 s_1, s_2 s_1 s_2\}$ we have $y \in \p^{m+1}$ and $x + \sqrt{d}y/2 \in \OF^\times + \P^{m+1}$. If $s' \in \{s_1, s_1 s_2, s_1 s_2 s_1, s_1 s_2 s_1 s_2\}$ we obtain $y \in \p^m$ and $x + \sqrt{d}y/2 \in \OF^\times + \P^m$. Hence, for any $s'$, we conclude that $g \in \GL_2(\OF)$. This gives us
\begin{align*}
& s' = 1, s_1 \Rightarrow X \in \mat{\p^{l+2m}}{\p^{l+m}}{\p^{l+m}}{\p^l}, \qquad s' = s_2 \Rightarrow X \in \mat{\p^{l+2m+1}}{\p^{l+m}}{\p^{l+m}}{\p^l},\\
&s' = s_2 s_1 \Rightarrow X \in \mat{\p^{l+2m+1}}{\p^{l+m+1}}{\p^{l+m+1}}{\p^l}, \qquad s' = s_1 s_2 \Rightarrow X \in \mat{\p^{l+2m}}{\p^{l+m}}{\p^{l+m}}{\p^{l+1}}, \\ 
& s' = s_1 s_2 s_1 \Rightarrow X \in \mat{\p^{l+2m}}{\p^{l+m+1}}{\p^{l+m+1}}{\p^{l+1}}, \qquad s' = s_2 s_1 s_2, s_1 s_2 s_1 s_2 \Rightarrow X \in \mat{\p^{l+2m+1}}{\p^{l+m+1}}{\p^{l+m+1}}{\p^{l+1}}.
\end{align*}
Now looking at the values of $B(h(l,m)s'), s' \in W$ from Proposition \ref{values-of-B-prop}, we get that either $B(s) = 0$ or $\Lambda(t) = \theta(u) = 1$.

We will illustrate one other case, $s = h(l,0) W_{w_0} s_1s_2$, since it is the most complicated. Here, $w_0$ is the unique element of $\OF/\p$ such that $w_0 + \alpha \not\in \OF_L^\times$. If $m_0 \geq 1$ or $l \leq -2$, then we have $B(s) = 0$. Hence, assume that $m_0=0$ and $l \geq -1$. Note that $x + y \frac{\sqrt{d}}2 = x-by/2-cw_0y+c(w_0 + \alpha)y$ and $a + bw_0 + cw_0^2 \in \p$. We see that $s^{-1} tu s \in \I$ implies that
$$y \in \OF, \quad x \pm (\frac b2 + cw_0)y \in \OF^\times.$$
 Hence, we see that $x + y \frac{\sqrt{d}}2 \in \OF_L^\times$. This implies that  $g \in \GL_2(\OF)$ and $\Lambda(t) = 1$. We have
$$\mat{1}{}{-w_0}{1}gX\mat{1}{-w_0}{}{1} \in \mat{\p^l}{\p^l}{\p^l}{\p^{l+1}}.$$
If $l \geq 0$, then we get $\theta(u) = 1$, as required. If $l = -1$, then let
$$\mat{1}{}{-w_0}{1}gX\mat{1}{-w_0}{}{1} = \mat{x_1}{x_2}{x_3}{x_4}, \mbox{ with } x_1, x_2, x_3 \in \varpi^{-1} \OF, x_4 \in \OF.$$
Set $\epsilon_1 = x + (b/2+cw_0) y, \epsilon_2 = x - (b/2+cw_0) y$. Using the fact that $X$ is symmetric and $\beta_{w_0}^0 \in \p$, we conclude that $x_3 \epsilon_1 - x_2 \epsilon_2 \in \OF$. Now $\theta(u) = \psi(\tr(SX))$ is equal to
\begin{align*}
& \psi(\frac 1{\det(g)} \Big(a((x-\frac{by}{2})x_1 - yc(x_3+w_0x_1)) + b(yax_1 + (x+\frac{by}{2})(x_3+w_0x_1)) \\
& \qquad + c(ya(x_2+w_0x_1) + (x+\frac{by}{2})(w_0^2x_1 + w_0(x_2+x_3) + x_4))\Big)) \\
&= \psi(\frac 1{\det(g)}\Big((x+\frac{by}{2})(x_1\beta_{w_0}^0 +c x_4) + x_2 \beta_{w_0}^0 yc - x_3 \beta_{w_0}^0 yc + (x_2 \epsilon_2 - x_3 \epsilon_1)c w_0 + x_3 \epsilon_1 (b+2c w_0) \Big)) \\
&= 1.
\end{align*}
Here, we have used that $x_3 \epsilon_1 - x_2 \epsilon_2 \in \OF, b+2cw_0 \in \p$ and $\psi$ is trivial on $\OF$. The other cases are computed in a similar manner (see \cite{P1}). \qed

\subsection{Criterion for $\dim(B(\Lambda, \theta)^\I)=1$}\label{criterion-for-dim-section}
In the previous sections, we have explicitly obtained a well-defined function $B$, which is right $\I$-invariant and satisfies (\ref{bessel-model-defining-condition}). The values of $B$ on the double coset representatives of $R(F) \backslash H(F) / \I$ were obtained, in Proposition \ref{values-of-B-prop}, using one or more of the conditions (\ref{hecke-operator-steinberg-1})-(\ref{hecke-operator-steinberg-2}). To show that the function $B$ is actually an element of $B(\Lambda, \theta)^\I$, we have to show that the conditions (\ref{hecke-operator-steinberg-1})-(\ref{hecke-operator-steinberg-2}) are satisfied by $B$ for every $h \in H(F)$. In fact, it is sufficient to show that $B$ satisfies  (\ref{hecke-operator-steinberg-1})-(\ref{hecke-operator-steinberg-2}) when $h$ is any double coset representative of $R(F) \backslash H(F) / \I$. The computations for checking this are long but not complicated (see \cite{P1}). We will describe the calculation for $h = h(l,m)$ below.

For $w, y \in \OF, w, y \neq 0$, we have the matrix identities
\begin{align}
\begin{bmatrix}1\\w&1\\&&1&-w\\&&&1\end{bmatrix} &= \begin{bmatrix}1&w^{-1}\\&1\\&&1\\&&-w^{-1}&1\end{bmatrix} s_1 \begin{bmatrix}-w\\&-w^{-1}\\&&-w^{-1}\\&&&-w\end{bmatrix} \begin{bmatrix}1&w^{-1}\\&1\\&&1\\&&-w^{-1}&1\end{bmatrix} \label{useful-matrix-identity-w} \\
\begin{bmatrix}1\\&1\\y&&1\\&&&1\end{bmatrix} &= \begin{bmatrix}1&&y^{-1}\\&1\\&&1\\&&&1\end{bmatrix} s_2 \begin{bmatrix}-y\\&1\\&&-y^{-1}\\&&&1\end{bmatrix} \begin{bmatrix}1&&y^{-1}\\&1\\&&1\\&&&1\end{bmatrix} \label{useful-matrix-identity-y}
\end{align}
Using (\ref{useful-matrix-identity-w}) and Lemmas \ref{coset-collapse-criterion-lemma}, \ref{beta-a-unit-lemma}, \ref{Lambda-character-sum-lemma}, we have
\begin{align*}
& \sum\limits_{w \in \OF/\p} B(h(l,m)s_1 W_w s_1) + B(h(l,m) s_1) = \sum\limits_{\substack{w \in \OF/\p \\ w \neq 0}} B(h(l,m) W_{w^{-1}} s_1) + B(h(l,m)) + B(h(l,m) s_1) \\
&= \sum\limits_{w \in \OF/\p} B(h(l,m) W_w s_1) + B(h(l,m)).
\end{align*}
By Proposition \ref{values-of-B-prop}, we see, for every value of $m_0, m, l, \Big(\frac{L}{\p}\Big)$, that the above quantity is equal to zero. Next, we have
\begin{align*}
B(h(l,m) \eta_0) &= B(h(l,m) \begin{bmatrix}\varpi\\&\varpi\\&&\varpi\\&&&\varpi\end{bmatrix}h(-1,0) s_2 s_1 s_2 ) \\
&= B(h(l-1,m) s_2 s_1 s_2) = \omega B(h(l,m)).
\end{align*}
Here, we have again used Proposition \ref{values-of-B-prop} and the identities $A_{l-1,m} = (-\omega q^3) A_{l,m}$. Finally, using (\ref{useful-matrix-identity-y}),
\begin{align*}
& \sum\limits_{y \in \OF/\p} B(h(l,m)\begin{bmatrix}1\\&1\\y&&1\\&&&1\end{bmatrix}) + B(h(l,m) s_2) \\
&= B(h(l,m)) + \sum\limits_{\substack{y \in \OF/\p \\ y \neq 0}} B(h(l,m) \begin{bmatrix}1&&y^{-1}\\&1\\&&1\\&&&1\end{bmatrix} s_2) + B(h(l,m) s_2) \\
&= B(h(l,m)) + \sum\limits_{y \in \OF/\p} \psi(a \varpi^{l+2m} y) B(h(l,m)s_2).
\end{align*}
By Proposition \ref{values-of-B-prop}, if $l \leq -1$, then both $B(h(l,m))$ and $B(h(l,m) s_2)$ are equal to zero, and if $l \geq 0$, then $B(h(l,m) s_2) = -1/q B(h(l,m))$. Hence, in either case, the above quantity is zero.

This shows that, for $h = h(l,m)$, the function $B$ satisfies (\ref{hecke-operator-steinberg-1}) - (\ref{hecke-operator-steinberg-2}), as required. The calculation for other values of $h$ follows in a similar manner (see \cite{P1}). Hence, we get the following theorem.
\begin{theorem}\label{B(lambda-theta)^I-dimension-theorem}
Let $\Lambda$ be a character of $L^\times$. Let $B(\Lambda, \theta)^\I$ be the space of smooth functions on $H(F)$, which are right $\I$-invariant, satisfy (\ref{bessel-model-defining-condition}) and the Hecke conditions (\ref{hecke-operator-steinberg-1}) - (\ref{hecke-operator-steinberg-2}). Then
\begin{equation}\label{dimension-equation}
{\rm dim}(B(\Lambda, \theta)^\I) = \left\{
                                     \begin{array}{ll}
                                       0, & \hbox{ if } \Lambda = \Omega \circ N_{L/F}  \mbox{ and } \Big(\frac L{\p}\Big) \in  \{-1, 0\}; \\
                                       1, & \hbox{ otherwise.}
                                     \end{array}
                                   \right.
\end{equation}
\end{theorem}
Note that the condition on $\Lambda$, in the case $\Big(\frac L{\p}\Big) \in  \{-1, 0\}$, follows from Proposition \ref{values-of-B-prop}, iv)a) and iv)c).

\subsection{Existence of Bessel model}\label{existence-section}
In this section we will obtain the existence of a $(\Lambda, \theta)$-Bessel model for $\pi$. In case $\Lambda$ is a unitary character, we will act with the Hecke algebra of $H(F)$ on a non-zero function in $B(\Lambda, \theta)^\I$. We will define an inner product on this Hecke module and also show that the Hecke module has a unique, up to a constant, function which is right $\I$-invariant (the same function that we started with). This will lead to the proof that the Hecke module is irreducible and is isomorphic to $\pi$, thus giving a $(\Lambda, \theta)$-Bessel model for $\pi$.

In case $\Lambda$ is not unitary (this can happen only if $L = F \oplus F$) we will obtain a Bessel model for $\pi$ using the Whittaker model.

\subsubsection*{Hecke module}
The Hecke algebra ${\mathcal H}$ of $H(F)$ is the space of all complex valued functions on $H(F)$ which are locally constant and compactly supported, with convolution product defined as follows,
\begin{equation}\label{convolution-product-defn}
(f_1 \ast f_2)(g) := \int\limits_{H(F)} f_1(h) f_2(h^{-1}g) dh, \qquad \mbox{ for } f_1, f_2 \in {\mathcal H}, g \in H(F).
\end{equation}
We refer the reader to \cite{Car} for details on Hecke algebras of $p$-adic groups and Hecke modules. Let $\Lambda$ be a character of $L^\times$ such that $B(\Lambda, \theta)^\I \neq 0$. Let $B \in B(\Lambda, \theta)^\I$ be the unique, up to a constant, function whose values are described in Proposition \ref{values-of-B-prop}. Define the action of $f \in {\mathcal H}$ on $B$ by
\begin{equation}\label{Hecke-algebra-action-on-B}
(R(f) B)(g) := \int\limits_{H(F)} f(h) B(gh) dh.
\end{equation}
This is a finite sum and hence converges for all $f$. Let
\begin{equation}\label{Hecke-module-defn}
V_B := \{R(f) B : f \in {\mathcal H} \}.
\end{equation}
Since $R(f_1) R(f_2) B = R(f_1 \ast f_2) B$, we see that $V_B$ is a Hecke module. Note that every function in $V_B$ transforms on the left according to $\Lambda \otimes \theta$.

\subsubsection*{Inner product on Hecke module}
Let us now assume that $\Lambda$ is a unitary character. Note that, by the comments in the begining of Sect.\ \ref{values-of-B-section}, if $L$ is a field, then $\Lambda$ is always unitary. In this case, we will define an inner product on the space $V_B$.
\begin{lemma}\label{norm-B-calculation-lemma}
The norm 
$$\langle B, B \rangle := \int\limits_{R(F) \backslash H(F)} |B(h)|^2 dh
$$ 
is finite. In particular,
\begin{equation}\label{norm-of-B-formula}
\langle B, B\rangle = \left\{
                        \begin{array}{ll}
                          {\rm vol}(\I) (1-\Big(\frac L{\p}\Big)q^{-1}) \frac{2q^{4m_0-3}}{(1-q^{-1})(1-q^{-3})} |C_{m_0}|^2, & \hbox{ if } m_0 \geq 1; \\
                          0 & \hbox{ if } m_0 = 0, \Big(\frac L{\p}\Big) = -1 \hbox{ OR } \\
                          & \quad \,\,m_0 = 0, \Big(\frac L{\p}\Big) = 0, \Lambda = \Omega \circ N_{L/F};\\
                          {\rm vol}(\I) \frac{2q^5+q^4+q^2-2q}{(1-q^{-3}) (1-q^{-1})} |C_0|^2, & \hbox{ if } m_0 = 0, \Big(\frac L{\p}\Big) = 0, \Lambda \neq \Omega \circ N_{L/F}; \\
                          {\rm vol}(\I) \frac{2(1+q^{-1})(q+2+q^{-1})}{1-q^{-3}} |B(W_{\frac{-b+\sqrt{d}}{2c}}s_1)|^2, & \hbox{ if } m_0 = 0, \Big(\frac L{\p}\Big) = 1, \omega \Lambda((1,\varpi)) = -1;\\
                          {\rm vol}(\I) \Big(\frac{2q^5}{1-q^{-3}} + \frac{2q^2(1-q^{-1})^3(q+2+q^{-1})}{(1-q^{-3})(1+\omega \Lambda((1,\varpi)))^2} \Big) |C_0|^2, & \hbox{ if } m_0 = 0, \Big(\frac L{\p}\Big) = 1, \omega \Lambda((1,\varpi)) \neq -1.
                        \end{array}
                      \right.
\end{equation}
Here, $C_m := B(h(0,m))$ and the measure is normalized so that ${\rm vol}(K^H) = 1$.
\end{lemma}
{\bf Proof.} We have
\begin{align*}
\int\limits_{R(F) \backslash H(F)} |B(h)|^2 dh &= \sum\limits_{s \in R(F) \backslash H(F) / \I} \,\,\,\int\limits_{R(F) \backslash R(F)s\I} |B(h)|^2 dh 
= \sum\limits_{s \in R(F) \backslash H(F) / \I}  |B(s)|^2 \int\limits_{\I_s \backslash \I} dh \\
&= \sum\limits_{s \in R(F) \backslash H(F) / \I} |B(s)|^2 \frac{{\rm vol}(\I)}{{\rm vol}(\I_s)}.
\end{align*}
Here, $\I_s := s^{-1}R(F)s \cap \I$. To get the last equality, we argue as in Lemma 3.7.1 of \cite{PS1}. The volume of $\I_s$ can be computed by similar methods to Sect.\ 3.7.1, 3.7.2 of \cite{PS1}. Now, using the values of $B(s)$ from Proposition \ref{values-of-B-prop} and geometric series, we get the result (see \cite{P1}). \qed

Let $L^2(R(F)\backslash H(F), \Lambda \otimes \theta) := \{ \phi : H(F) \rightarrow \C : \mbox{ smooth, }  \phi(rh) = \Lambda \otimes \theta(r) \phi(h) \mbox{ for } r \in R(F), h \in H(F), \int_{R(F)\backslash H(F)} \mid \phi(h)\mid^2 dh < \infty \}$. The previous lemma tells us that $B \in L^2(R(F)\backslash H(F), \Lambda \otimes \theta)$. It is an easy exercise to see that, in fact, for any $f \in \mathcal H$, we have $R(f)B \in L^2(R(F)\backslash H(F), \Lambda \otimes \theta)$. Now, we see that $V_B$ inherits the inner product from $L^2(R(F)\backslash H(F), \Lambda \otimes \theta)$. For $f_1, f_2 \in {\mathcal H}$, we obtain
\begin{equation}\label{norm-on-V-defn}
\langle R(f_1)B, R(f_2)B \rangle = \int\limits_{R(F) \backslash H(F)} (R(f_1)B)(g) \overline{(R(f_2)B)(g)} dg.
\end{equation}

\begin{lemma}\label{adjoint-lemma}
For $f \in \mathcal H$, define $f^\ast \in \mathcal H$ by $f^\ast(g) = \overline{f(g^{-1})}$. Then we have
\begin{equation}\label{adjoint-eqn}
\langle B_1, R(f)B_2 \rangle = \langle R(f^\ast) B_1, B_2 \rangle, \qquad \mbox{ for any } B_1, B_2 \in V_B.
\end{equation}
\end{lemma}
{\bf Proof.} The lemma follows by a formal calculation. \qed

\subsection*{Irreducibility of $V_B$}
\begin{lemma}\label{V-B^I-is-one-dimensional-lemma}
Let $V_B^\I$ be the subspace of functions in $V_B$ that are right $\I$-invariant.
Then
$${\rm dim}(V_B^\I) = 1.$$
\end{lemma}
{\bf Proof.} We know that $V_B^\I$ is not trivial since $B \in V_B^\I$. Let $\chi_\I \in {\mathcal H}$ be the characteristic function of $\I$ and set $f_\I := {\rm vol}(\I)^{-1} \chi_\I$. Then, by definition, any $B' \in V_B^I$, satisfies $R(f_\I) B' = B'$. Let $f \in \mathcal H$ be such that $B' = R(f) B = R(f \ast f_\I) B$. Here, we have used that $B \in V_B^\I$. Then
$$B' = R(f_\I) B' = R(f_\I) (R(f \ast f_\I) B) = R(f_\I \ast f \ast f_\I) B.$$
But $f_\I \ast f \ast f_\I \in \mathcal H_\I$, the Iwahori Hecke algebra. Since $B$ is an eigenfunction of $\mathcal H_\I$, we see that $B' \in \C B$. Hence, ${\rm dim}(V_B^\I) = 1$, as required. \qed

\begin{proposition}\label{irreducibility-of-VB}
Let $\pi = \Omega {\rm St}_{\GSp_4}$ be the Steinberg representation of $H(F)$, twisted by an unramified, quadratic character $\Omega$. Let $\Lambda$ be a character of $L^\times$ such that $\dim(B(\Lambda, \theta)^\I) = 1$. Let $V_B$ be as in (\ref{Hecke-module-defn}). If $\Lambda$ is unitary, then $V_B$ is irreducible and isomorphic to $\pi$.
\end{proposition}
{\bf Proof.} Let us assume, to the contrary, that $V_B$ is reducible. Let $W$ be an $\mathcal H$-invariant subspace. Let $W^{\perp}$ be the complement of $W$ in $V_B$ with respect to the inner product $\langle \, , \, \rangle$ defined in (\ref{norm-on-V-defn}). Using Lemma \ref{adjoint-lemma}, we see that $W^\perp$ is also $\mathcal H$-invariant. Write $B = B_1 + B_2$, with $B_1 \in W, B_2 \in W^\perp$. Let $f_\I$ be as defined in the proof of Lemma \ref{V-B^I-is-one-dimensional-lemma}. Since $W, W^\perp$ are $\mathcal H$-invariant, we see that $R(f_\I) B_1 \in W$ and $R(f_\I) B_2 \in W^\perp$.  Since $B$ is right $\I$-invariant, we see that $B_1 = R(f_\I) B_1$ and $B_2 = R(f_\I) B_2$. By Lemma \ref{V-B^I-is-one-dimensional-lemma}, we obtain, either $B = B_1$ or $B = B_2$. Since $V_B$ is generated by $B$, we have either $W = V_B$ or $W = 0$. Hence, we see that $V_B$ is an irreducible Hecke module, which contains a unique, up to a constant, vector which is right $\I$-invariant. This uniquely characterizes the Steinberg representation of $H(F)$, and hence, $V_B$ is isomorphic to $\pi$. \qed

\subsubsection*{Generic representations have split Bessel models}
Let us now assume that $\Lambda$ is not a unitary character. This can happen only if $L = F \oplus F$. In this case, we will use the fact that $\Omega {\rm St}_{\GSp_4}$ is a generic representation. We will now show that any irreducible, admissible, generic representation of $H(F)$ has a split Bessel model. 

Let $S = \mat{a}{b/2}{b/2}{c}$ be such that $b^2-4ac$ is a square in $F^\times$. One can find a matrix $A \in \GL_2(\OF)$ such that $S' :=  \,^tA S A = \mat{}{1/2}{1/2}{}$. In this case, $T_{S'}(F) := \{g \in \GL_2(F) : \,^tg S' g = \det(g) S'\} = A^{-1} T(F) A.$ The group $T_{S'}(F)$ embedded in $H(F)$ is given by
$$\{\begin{bmatrix}x\\&y\\&&y\\&&&x\end{bmatrix} : x, y \in F^\times \} .$$
Let $\theta'$ be the character of $U(F)$ obtained from $S'$ and $\Lambda'$ be the character of $T_{S'}(F)$ obtained from $\Lambda$. Then it is easy to see that $\pi$ has a $(\Lambda, \theta)$-Bessel model if and only if it has a $(\Lambda', \theta')$-Bessel model. So, we will assume that $S = \mat{}{1/2}{1/2}{}$.

Let $(\pi, V)$ be an irreducible, admissible representation of $H(F)$.  For $c_1,c_2\in F^\times$, consider the character $\psi_{c_1,c_2}$ of the
unipotent radical $N(F)$ of the Borel subgroup given by
$$
 \psi_{c_1,c_2}(\begin{bmatrix}1&x&*&*\\&1&*&y\\&&1&\\&&-x&1\end{bmatrix})
 =\psi(c_1x+c_2y).
$$
The representation $\pi$ of $H(F)$ is called \emph{generic}
if ${\rm Hom}_{N(F)}(\pi,\psi_{c_1,c_2})\neq0$. In this case there is an
associated Whittaker model $\mathcal{W}(\pi,\psi_{c_1,c_2})$ consisting of functions
$H(F)\rightarrow\C$ that transform on the left according to $\psi_{c_1,c_2}$.
For $W\in\mathcal{W}(\pi,\psi_{c_1,c_2})$, there is an associated zeta integral
\begin{equation}\label{ZsWdefeq}
 Z(s,W)=\int\limits_{F^\times}\int\limits_F
 W(\begin{bmatrix}y\\&y\\&&1\\&x&&1\end{bmatrix})|y|^{s-3/2}\,dx\,d^\times y.
\end{equation}
This integral is convergent for ${\rm Re}(s)>s_0$, where $s_0$ is independent
of $W$ (\cite{RS}, Proposition 2.6.3). More precisely, the integral converges
to an element of $\C(q^{-s})$, and therefore has meromorphic continuation
to all of $\C$. Moreover, there exists an $L$-factor of the form
$$
 L(s,\pi)=\frac1{Q(q^{-s})},\qquad Q(X)\in\C[X],\;Q(0)=1,
$$
such that
\begin{equation}\label{zetaLquotienteq}
 \frac{Z(s,W)}{L(s,\pi)}\in\C[q^{-s},q^s]\qquad\text{for all }
 W\in \mathcal{W}(\pi,\psi_{c_1,c_2}).
\end{equation}
(This is proved in \cite{RS} Proposition 2.6.4 for $\pi$ with trivial central
character.)
\begin{lemma}\label{GSp4genericlemma}
 Let $(\pi,V)$ be an irreducible, admissible, generic representation of $H(F)$ with trivial central character.
 Let $\sigma$ be a unitary character of $F^\times$, and let $s\in\C$ be arbitrary.
 Then there exists a non-zero functional $f_{s,\sigma}:\:V\rightarrow\C$
 with the following properties.
 \begin{enumerate}
  \item For all $x,y,z\in F$ and $v\in V$,
   \begin{equation}\label{GSp4genericlemmaeq1}
    f_{s,\sigma}(\pi(\begin{bmatrix}1&&x&y\\&1&y&z\\&&1\\&&&1\end{bmatrix})v)
    =\psi(c_1y)f_{s,\sigma}(v).
   \end{equation}
  \item For all $x\in F^\times$ and $v\in V$,
   \begin{equation}\label{GSp4genericlemmaeq2}
    f_{s,\sigma}(\pi(\begin{bmatrix}x\\&1\\&&1\\&&&x\end{bmatrix})v)
    =\sigma(x)^{-1}|x|^{-s+1/2}f_{s,\sigma}(v).
   \end{equation}
 \end{enumerate}
\end{lemma}
{\bf Proof:} We may assume that $V=\mathcal{W}(\pi,\psi_{c_1,c_2})$. Let $s_0\in\R$
be such that $Z(s,W)$ is absolutely convergent for ${\rm Re}(s)>s_0$.
Then the integral
\begin{equation}\label{ZsWsigmadefeq}
 Z_\sigma(s,W)=\int\limits_{F^\times}\int\limits_F
 W(\begin{bmatrix}y\\&y\\&&1\\&x&&1\end{bmatrix})|y|^{s-3/2}\sigma(y)\,dx\,d^\times y
\end{equation}
is also absolutely convergent for ${\rm Re}(s)>s_0$, since $\sigma$ is unitary.
Note that these are the zeta integrals for the twisted representation $\sigma\pi$.
Therefore, by (\ref{zetaLquotienteq}), the quotient $Z_\sigma(s,W)/L(s,\sigma\pi)$
is in $\C[q^{-s},q^s]$ for all $W\in\mathcal{W}(\pi,\psi_{c_1,c_2})$.
Now, for ${\rm Re}(s)>s_0$, we define
\begin{equation}\label{GSp4genericlemmaeq3}
 f_{s,\sigma}(W)=\frac{Z_\sigma(s,\pi(w)W)}{L(s,\sigma\pi)}, \qquad \mbox{ where } w=\begin{bmatrix}1\\&&&1\\&&1\\&-1&&\end{bmatrix}.
\end{equation}
Straightforward calculations show that (\ref{GSp4genericlemmaeq1}) and
(\ref{GSp4genericlemmaeq2}) are satisfied. For general $s$, since the
quotient (\ref{GSp4genericlemmaeq3}) is entire, we can define $f_{s,\sigma}$ by
analytic continuation.\qed

\begin{proposition}\label{GSp4genericprop}
 Let $(\pi,V)$ be an irreducible, admissible, generic representation of $H(F)$ with trivial central character. Then $\pi$ admits a split Bessel functional with respect to any character
 $\Lambda$ of $T(F)$ that satisfies $\Lambda\big|_{F^\times} \equiv 1$.
\end{proposition}
{\bf Proof:} As mentioned earlier, we can take $S = \mat{}{1/2}{1/2}{}$.
Let $s\in\C$ and $\sigma$ be a unitary character of $F^\times$ such that
$$
 \Lambda(\begin{bmatrix}x\\&1\\&&1\\&&&x\end{bmatrix})=\sigma(x)^{-1}|x|^{-s+1/2}
 \qquad\text{for all }x\in F^\times.
$$
Let $f_{s,\sigma}$ be as in Lemma \ref{GSp4genericlemma}. We may assume that
$c_1=1$, so that $f_{s,\sigma}(\pi(u)v)=\theta(u)v$ for all $u\in U(F)$ by
(\ref{GSp4genericlemmaeq1}). We have
\begin{equation}\label{GSp4genericpropeq1}
 f_{s,\sigma}(\pi(\begin{bmatrix}x\\&1\\&&1\\&&&x\end{bmatrix})v)
 =\Lambda(x)f_{s,\sigma}(v)\qquad\text{for all }x\in F^\times
\end{equation}
by (\ref{GSp4genericlemmaeq2}). Since $\Lambda\big|_{F^\times} \equiv 1$ we
in fact obtain $f_{s,\sigma}(\pi(t)v)=\Lambda(t)f_{s,\sigma}(v)$ for all $t\in T(F)$. Hence
$f_{s,\sigma}$ is a Bessel functional as desired.\qed

Let us remark here that, in the split case, for values of $s \in \C$ outside the range of convergence of the zeta integral, we do not have an explicit formula for the Bessel functional. This, in turn, is also reflected in the fact that it is not very easy to define an inner product on the space $V_B$ (defined in (\ref{Hecke-module-defn})), although it is known that the Steinberg representation is square-integrable.

\subsubsection*{Main result on existence and uniqueness of Bessel models}
\begin{theorem}\label{existence-theorem}
Let $\pi = \Omega {\rm St}_{\GSp_4}$ be the Steinberg representation of $H(F)$, twisted by an unramified quadratic character $\Omega$. Let $\Lambda$ be a character of $L^\times$ such that $\Lambda \mid_{F^\times} \equiv 1$. If $L$ is a field, then $\pi$ has a $(\Lambda, \theta)$-Bessel model if and only if $\Lambda \neq \Omega \circ N_{L/F}$. If $L$ is not a field, then $\pi$ always has a $(\Lambda, \theta)$-Bessel model. In case $\pi$ has a $(\Lambda, \theta)$-Bessel model, it is unique. 

In addition, if $\pi$ has a $(\Lambda, \theta)$-Bessel model, then the Iwahori spherical vector of $\pi$ is a test vector for the Bessel functional if and only if $\Lambda$ satisfies the following conditions.
\begin{enumerate}
\item $\Lambda \mid_{1+\P} \equiv 1$, i.e., $c(\Lambda) \leq 1$ (see (\ref{c(Lambda)-defn}) for definition of $c(\Lambda)$).
\item If $\Big(\frac L\p\Big) = 1$ and $\Lambda$ is unramified, then $\Lambda((1,\varpi)) \neq \Omega(\varpi)$.
\end{enumerate}
\end{theorem}
{\bf Proof.} If $\pi$ has a $(\Lambda, \theta)$-Bessel model, then it contains a unique vector in $B(\Lambda, \theta)^\I$. By Theorem \ref{B(lambda-theta)^I-dimension-theorem}, the dimension of $B(\Lambda, \theta)^\I$ is one, which gives us the uniqueness of Bessel models.

Now we will show the existence of the Bessel model. Let $\Lambda$ be a character of $L^\times$, with $\Lambda \mid_{F^\times} \equiv 1$, such that, if $L$ is a field, $\Lambda \neq \Omega \circ N_{L/F}$ . We know, by Proposition \ref{B(lambda-theta)^I-dimension-theorem}, that ${\rm dim}(B(\Lambda, \theta)^\I) = 1$.  If $\Lambda$ is unitary, then Proposition \ref{irreducibility-of-VB} tells us that $V_B$ is a $(\Lambda, \theta)$-Bessel model for $\pi$. If $\Lambda$ is not unitary, then we use the fact that $\pi$ is a generic representation in the split case. Then Proposition \ref{GSp4genericprop} gives us the result. 

The statement regarding the test vector can be deduced from Proposition \ref{values-of-B-prop} and the fact that a Bessel function $B$ corresponds to a test vector if and only if $B(1) \neq 0$.  \qed

\section{Integral representation of the non-archimedean local $L$-function}\label{non-arch-int-rep-section}
In this section, using the explicit values of the Bessel function obtained in Proposition \ref{values-of-B-prop}, we will obtain an integral representation of the $L$-function for the Steinberg representation $\pi$ of $H(F)$ twisted by any irreducible, admissible representation $\tau$ of $\GL_2(F)$. For this, we will use the integral obtained by Furusawa in \cite{Fu}. Let us briefly describe the setup.
\subsection{The unitary group, parabolic induction and the local integral}\label{local-induced-repn-section}
Let $G = \GU(2,2;L)$ be the unitary similitude group, whose $F$-points are given by
\begin{equation}\label{unitary-gp-defn}
G(F) := \{g \in \GL_4(L) : \, ^t\bar{g}Jg = \mu_2(g)J,\:\mu_2(g)\in F^{\times}\},
\end{equation}
where $J = \mat{}{1_2}{-1_2}{}$. Note that $H(F) = G(F)
\cap\GL_4(F)$. As a minimal parabolic subgroup we choose the
subgroup of all matrices that become upper triangular after
switching the last two rows and last two columns. Let $P$ be the
standard maximal parabolic subgroup of $G(F)$ with a non-abelian
unipotent radical. Let $P = MN$ be the Levi decomposition of $P$.
We have $M = M^{(1)}M^{(2)}$, where
\begin{eqnarray}
&& M^{(1)}(F)=\{\begin{bmatrix}\zeta\\&1\\&&\bar
  \zeta^{-1}\\&&&1\end{bmatrix}:\:\zeta\in L^\times\}, \quad
 M^{(2)}(F)=\{\begin{bmatrix}1\\&\alpha&&\beta\\&&\mu\\&\gamma&&\delta\end{bmatrix}
  \in G(F)\} \nonumber,\\
&& N(F) =\{\begin{bmatrix}
                 1 & z &  &  \\
                  & 1 &  &  \\
                  &  & 1 &  \\
                  &  & -\overline{z} & 1 \\
               \end{bmatrix}
               \begin{bmatrix}
                 1 &  & w & y \\
                  & 1 & \overline{y} &  \\
                  &  & 1 &  \\
                  &  &  & 1 \\
               \end{bmatrix}
              : w\in F,\;y,z \in L \}. \label{M1-M2-N-defn}
\end{eqnarray}
The modular factor of the parabolic $P$ is given by
\begin{equation}\label{deltaPformulaeq}
 \delta_P(\begin{bmatrix}\zeta\\&1\\&&\bar\zeta^{-1}\\&&&1\end{bmatrix}
 \begin{bmatrix}1\\&\alpha&&\beta\\&&\mu\\&\gamma&&\delta\end{bmatrix})
 =|N(\zeta)\mu^{-1}|^3\qquad(\mu=\bar\alpha\delta-\beta\bar\gamma),
\end{equation}
where $|\cdot|$ is the normalized absolute value on $F$.
Let $(\tau,V_\tau)$ be an irreducible, admissible representation
of $\GL_2(F)$, and let $\chi_0$ be a character of $L^\times$ such
that $\chi_0\big|_{F^\times}$ coincides with $\omega_{\tau}$, the
central character of $\tau$. Let us assume that $V_\tau$ is the Whittaker model of $\tau$ with respect to the character $\psi^{-c}$ (we assume that $c \neq 0$). Then the representation
$(\lambda,g)\mapsto\chi_0(\lambda)\tau(g)$ of
$L^\times\times\GL_2(F)$ factors through
$\{(\lambda,\lambda^{-1}):\:\lambda\in F^\times\}$, and
consequently defines a
representation of $M^{(2)}(F)$ on the same space $V_\tau$. Let $\chi$ be a character of $L^\times$, considered as a character of $M^{(1)}(F)$. Extend the representation $\chi \times \chi_0 \times \tau$ of $M(F)$ to a representation of $P(F)$ by setting it to be trivial on $N(F)$. If $s$ is a complex parameter, set $I(s,\chi,\chi_0,\tau)= {\rm Ind}_{P(F)}^{G(F)}(\delta_P^{s+1/2} \times \chi \times \chi_0 \times \tau)$.

Let $(\pi,V_\pi)$ be the twisted Steinberg representation of
$H(F)$. We assume that $V_\pi$ is a Bessel model for $\pi$
with respect to a character $\Lambda\otimes\theta$ of $R(F)$.
Let the characters $\chi, \chi_0$ and $\Lambda$ be related by
\begin{equation}\label{chi-lambda-char-condition}
\chi(\zeta) = \Lambda(\bar{\zeta})^{-1} \chi_0(\bar{\zeta})^{-1}.
\end{equation}
 Let
$W^\#(\,\cdot\,,s)$ be an element of $I(s,\chi,\chi_0,\tau)$ for
which the restriction of $W^\#(\,\cdot\,,s)$ to the standard
maximal compact subgroup of $G(F)$  is
independent of $s$, i.e., $W^\#(\,\cdot\,,s)$ is a ``flat section'' of
the family of induced representations $I(s,\chi,\chi_0,\tau)$. By
Lemma 2.3.1 of \cite{PS1}, it is meaningful to consider the
integral
\begin{equation}\label{localZseq}
 Z(s)=\int\limits_{R(F)\backslash H(F)}W^\#(\eta h,s)B(h)\,dh,
\end{equation}
where
\begin{equation}\label{eta-defn}
\eta = \begin{bmatrix}1&&&\\\alpha&1&&\\&&1&-\bar\alpha\\&&&1\end{bmatrix}.
\end{equation}
This is the local component of the global integral considered in Sect.\ \ref{global-repn-integral} below.

\subsection{The $\GL_2$ newform}\label{the-GL2-newform-section}
 Let us define $K^{(0)}(\p^0) = \GL_2(\OF)$ and,
for $n>0$,
\begin{equation}\label{K1defeq}
 K^{(0)}(\p^n)=\GL_2(\OF)\cap\mat{1+\p^n}{\OF}{\p^n}{\OF^\times}.
\end{equation}
As above, let $(\tau,V_\tau)$ be a generic, irreducible, admissible
representation of $\GL_2(F)$ such that $V_\tau$ is the
$\psi^{-c}$--Whittaker model of $\tau$. It is well known that $V_\tau$ has a unique (up to a constant) vector $W^{(1)}$, called the newform, that is right-invariant under $K^{(0)}(\p^n)$ for some $n \geq 0$. We then say that $\tau$ has conductor $\p^n$. Let us normalize $W^{(1)}$ so that $W^{(1)}(1) = 1$. We will need the values of $W^{(1)}$ evaluated at $\mat{\varpi^l}{}{}{1}$, for $l \geq 0$. The following table gives these values (refer Sect.\ 2.4 \cite{Sc1}). 

$$\renewcommand{\arraystretch}{1.4}
 \begin{array}{|c|c|}
 \hline 
\tau &  W^{(1)}(\mat{\varpi^l}{}{}{1})\\ \hline 
\alpha \times \beta \mbox{ with } \alpha, \beta \mbox{ unramified, } \alpha \beta^{-1} \neq \mid \, \mid^{\pm 1} & q^{-\frac l2} \frac{\alpha(\varpi^{l+1}) - \beta(\varpi^{l+1})}{\alpha(\varpi) - \beta(\varpi)} \\ \hline
\alpha \times \beta \mbox{ with } \alpha \mbox{ unramified, } \beta \mbox{ ramified, } \alpha \beta^{-1} \neq \mid \, \mid^{\pm 1} & \omega_\tau(\varpi^l)
    \alpha(\varpi^{-l}) q^{-\frac l2} \\ \hline
\mbox{ supercuspidal OR ramified twist of Steinberg}  & 1 \qquad \mbox{ if } l = 0 \\
\mbox{ OR } \alpha \times \beta \mbox{ with } \alpha, \beta \mbox{ ramified, } \alpha \beta^{-1} \neq \mid \, \mid^{\pm 1}  & 0 \qquad \mbox{ if } l > 0 \\ \hline
\Omega' {\rm St}_{\GL_2}, \mbox{ with } \Omega' \mbox{ unramified } & \Omega'(\varpi^l)q^{-l} \\ \hline
\end{array}
$$

We extend $W^{(1)}$ to a function on $M^{(2)}(F)$ via
\begin{equation}\label{W1extensioneq}
 W^{(1)}(ag)=\chi_0(a)W^{(1)}(g),\qquad a\in L^\times,\:g\in\GL_2(F).
\end{equation}

\subsection{Choice of $\Lambda$ and $W^\#$}\label{choice-of-W}
We will choose a character $\Lambda$ of $L^\times$ such that $\pi$ has a $(\Lambda, \theta)$-Bessel model and the Iwahori spherical vector is a test vector for the Bessel functional. Noting that $\Lambda \mid_{F^\times}$ is the central character of $\pi$ and using Theorem \ref{existence-theorem}, we impose the following conditions on $\Lambda$.
\begin{enumerate}
\item  $\Lambda |_{F^\times} \equiv 1$

\item If $L$ is a field, then $\Lambda \neq \Omega \circ N_{L/F}$

\item $c(\Lambda) \leq 1$ 

\item If $L$ is not a field and $c(\Lambda) = 0$, then $\omega \Lambda((1,\varpi)) \neq -1$.
\end{enumerate}
Note that this implies that $\Lambda |_{\OF^\times + \P} \equiv 1$. For $n \geq 1$, let $\Gamma(\P^n)$ be the principal congruence subgroup of the maximal compact subgroup $K^G := G(\OF)$ of $G(F)$, defined by
\begin{equation}\label{princ-cong-defn-steinberg}
 \Gamma(\P^n) := \{g \in K^G : g \equiv 1 \pmod{\P^n}\}.
\end{equation}
We prove the following lemma, which will be crucial for the well-definedness of $W^\#$ below.
\begin{lemma}\label{well-defined-sub-lemma}
 Let $(\tau,V_\tau)$ be a generic, irreducible, admissible
representation of $\GL_2(F)$ with conductor $\p^n, n \geq 0$. Set $n_0 = {\rm max}\{1, n\}$ and let
 $$
  \hat{m} = \begin{bmatrix}\zeta\\&a'&&b'\\&&\mu\bar{\zeta}^{-1}\\&c'&&d'\end{bmatrix} \in M(F)
  \qquad\mbox{and}\qquad \hat{n} = \begin{bmatrix}
                 1 & z &  &  \\
                  & 1 &  &  \\
                  &  & 1 &  \\
                  &  & -\overline{z} & 1 \\
               \end{bmatrix}
               \begin{bmatrix}
                 1 &  & w & y \\
                  & 1 & \overline{y} &  \\
                  &  & 1 &  \\
                  &  &  & 1 \\
               \end{bmatrix} \in N(F).
 $$
 Suppose we have $A := \eta^{-1} \hat{m}\hat{n} \eta \in \I \Gamma(\P^{n_0})$.  Then  we get
\begin{enumerate}
\item $c' \in \P^{n_0}$ and $a' \bar{\zeta}^{-1}  \in 1+\P^{n_0}$.

\item for any  $\mat{a_1'}{b_1'}{c_1'}{d_1'} \in \GU(1,1;L)(F)$,
$$\chi(\zeta) W^{(1)}(\mat{a_1'}{b_1'}{c_1'}{d_1'} \mat{a'}{b'}{c'}{d'}) = W^{(1)}(\mat{a_1'}{b_1'}{c_1'}{d_1'}).$$
\end{enumerate}
\end{lemma}

{\bf Proof.} Using Lemma \ref{quadraticextdisclemma-steinberg} i), it is easy to show that for $n \geq 0$
\begin{equation}\label{something-in-pn}
x \in \OF+\P^n \mbox{ and } \alpha x \in \OF+\P^n \mbox{ implies } x \in \P^n.
\end{equation}
First note that $\I \Gamma(\P^{n_0}) \subset
M_4(\OF+\P^{n_0})$. Looking at the $(4,1), (4,2)$ coefficient
of $A$, we see that $c', \alpha c' \in
\OF+\P^{n_0}$. By (\ref{something-in-pn}), we obtain $c'
\in \P^{n_0}$, as required.

Observe that $\hat{m} \hat{n} \in K^G$ and $c' \in \P^{n_0} \subset \P$ implies that $\zeta, a', d' \in \OF_L^\times$. The upper left $2 \times 2$
block of $A$ is given by
$$\mat{\zeta + \alpha z \zeta}{z\zeta}{\alpha a' -\alpha (\zeta + \alpha z \zeta ))}{a'-\alpha z
\zeta}.$$
We will repeatedly use the following fact:
\begin{equation}\label{well-defined-sub-lemmaeq1}
 \text{If }x \in \OF + \P^{n_0},\text{ then }x \equiv \bar{x} \pmod{(\alpha-\bar\alpha)\P^{n_0}}.
\end{equation}
For, if $x=y+\alpha z$ with $y\in\OF$ and $z\in\p^{n_0}$, then $x-\bar x=(\alpha-\bar\alpha)z$.
Applying this to the matrix entries of $A$, we get $z\zeta \equiv \bar{z} \bar{\zeta} \pmod{(\alpha-\bar\alpha)\P^{n_0}}$, and then
\begin{equation}\label{a'-zeta-condition}
 a'-\bar{a'} \equiv (\alpha - \bar{\alpha}) z \zeta
 \pmod{(\alpha-\bar\alpha)\P^{n_0}}, \qquad \zeta - \bar{\zeta} \equiv
 (\bar{\alpha}-\alpha) z \zeta  \pmod{(\alpha-\bar\alpha)\P^{n_0}}.
\end{equation}
Using $\zeta + \alpha z \zeta \equiv \bar\zeta + \bar\alpha \bar{z} \bar\zeta
 \pmod{(\alpha-\bar\alpha)\P^{n_0}}$ and (\ref{a'-zeta-condition}),
we get from the $(2,1)$ coefficient of $A$ that
$$
 (a' - \bar\zeta) (\alpha - \bar\alpha) \equiv 0
  \pmod{(\alpha-\bar\alpha)\P^{n_0}}.
$$
Hence $a' - \bar\zeta \equiv0\pmod{\P^{n_0}}$,
so that $a' \bar{\zeta}^{-1} \in 1+\P^{n_0}$, as required. This proves part i) of the lemma.

Looking at the $(1,2)$ coefficient of $A$, we see that $z \zeta \in \P$. Looking at the $(1,1)$ coefficient of $A$, we see that $\zeta \in \OF^\times + \P$.
\begin{eqnarray*}
\chi(\zeta) W^{(1)}(\mat{a_1'}{b_1'}{c_1'}{d_1'} \mat{a'}{b'}{c'}{d'}) &=& \chi(\zeta) \chi_0(a') W^{(1)}(\mat{a_1'}{b_1'}{c_1'}{d_1'} \mat{1}{b'/a'}{c'/a'}{d'/a'}) \\
&=&  \Lambda(\bar{\zeta}^{-1}) \chi_0(\bar{\zeta}^{-1}) \chi_0(a') W^{(1)}(\mat{a_1'}{b_1'}{c_1'}{d_1'} \mat{1}{b'/a'}{c'/a'}{d'/a'}) \\
&=& W^{(1)}(\mat{a_1'}{b_1'}{c_1'}{d_1'})
\end{eqnarray*}
Here, we have used the fact that $\Lambda$ is trivial on $\OF^\times + \P, \chi_0$ is trivial on $1 + \P^{n_0}$ and the matrix $\mat{1}{b'/a'}{c'/a'}{d'/a'}$ lies in $K^{(0)}(\p^{n_0})$. \qed

Let $n_0 = {\rm max}\{1, n\}$, as above. Given a complex number $s$, define the function
$W^\#(\,\cdot\,,s):\:G(F)\rightarrow\C$ as follows.
\begin{enumerate}
 \item If $g\notin M(F)N(F)\eta {\rm I} \Gamma(\P^{n_0})$, then $W^\#(g,s)=0$.
 \item If $g=mn\eta k \gamma$ with $m\in M(F)$, $n\in N(F)$, $k\in {\rm I}$,
  $\gamma \in \Gamma(\P^{n_0})$, then $W^\#(g,s)=W^\#(m \eta,s)$.
 \item For $\zeta\in L^\times$ and $\mat{a'}{b'}{c'}{d'}\in M^{(2)}(F)$,
  \begin{equation}\label{Wsharpformulaeq-steinberg}
   W^\#(\begin{bmatrix}\zeta\\&1\\&&\bar{\zeta}^{-1}\\&&&1\end{bmatrix}
   \begin{bmatrix}1\\&a'&&b'\\&&\mu\\&c'&&d'\end{bmatrix} \eta,s)
   =|N(\zeta)\cdot\mu^{-1}|^{3(s+1/2)}\chi(\zeta)\,
   W^{(1)}(\mat{a'}{b'}{c'}{d'}).
  \end{equation}
  Here $\mu=\bar{a'}d'-b'\bar{c'}$.
\end{enumerate}
By Lemma \ref{well-defined-sub-lemma}, we see that $W^\#$ is well-defined. It is an element of $I(s,\chi,\chi_0,\tau)$.

\subsection{Support of $W^\#$}
Let us choose $W^\#$ as above and $B$ as in Proposition \ref{values-of-B-prop}, with $B(1) = 1$. Note that $B(1) \neq 0$ by the comments in the begining of Sect.\ \ref{choice-of-W}. Then the integral (\ref{localZseq}) becomes
\begin{equation}\label{integral-as-a-sum-1}
Z(s) = \sum\limits_{l \in \Z, m \geq 0} \sum\limits_{t} W^\#(\eta h(l,m) t, s) B(h(l,m)t) V^{l,m}_t,
\end{equation}
where $t$ runs through the double coset representatives from Proposition \ref{disjoint-double-cosets-prop} and
$$V_t^{l,m} = {\rm vol}(R(F) \backslash R(F) h(l,m) t \I).$$
To compute (\ref{integral-as-a-sum-1}), we need to find out for what values of $l, m, t$ is $\eta h(l,m) t$ in the support of $W^\#$. Write $\eta h(l,m) = h(l,m) \eta_m$, where
\begin{equation}\label{etam-defn}
\eta_m = \begin{bmatrix}1&&&\\ \varpi^m\alpha&1&&\\&&1&-\varpi^m\bar\alpha\\&&&1\end{bmatrix}.
\end{equation}
Since $h(l,m) \in M(F)$, we need to know for which values of $m, t$ is $\eta_m t$ in the support of $W^\#$. This is done in the following lemma.
\begin{lemma}\label{etam-t-in-supp-Wsharp}
Let $t$ be any double coset representative from Proposition \ref{disjoint-double-cosets-prop}. Then $\eta_m t$ lies in the support, $M N \eta {\rm I} \Gamma(\P^{n_0})$, of $W^\#$ if and only if $m=0$ and $t=1$.
\end{lemma}
{\bf Proof.} Let us first consider the case $m > 0$.  Note that it is enough to show that $\eta_m t \notin M N \eta {\rm I} \Gamma(\P)$. For any double coset representative $t$, we have $t^{-1} \eta_m t \equiv 1 \pmod{\P}$ and hence $t^{-1} \eta_m t \in \Gamma(\P)$. So it is enough to show that $t \notin M N \eta {\rm I} \Gamma(\P)$ for any $t$. Suppose, there are $\hat{m} \in M, \hat{n} \in N$ such that $A = \eta^{-1} \hat{m} \hat{n} t \in {\rm I} \Gamma(\P)$. Using $\hat{m}, \hat{n} \in K^G$ and
\begin{equation}\label{I-Gamma-description}
{\rm I} \Gamma(\P) \subset \begin{bmatrix}\OF+\P&\P&\OF+\P&\OF+\P\\\OF+\P&\OF+\P&\OF+\P&\OF+\P\\\P&\P&\OF+\P&\OF+\P\\
\P&\P&\P&\OF+\P\end{bmatrix}
\end{equation}
we get a contradiction for every $t \in W$. Let us now consider the case $m = 0$. First let $t=1$. Taking $\hat{m} = \hat{n} = 1$, we easily see that $\eta \in M N \eta {\rm I} \Gamma(\P^{n_0})$, as required. Now assume that $t \neq 1$. Suppose, there are $\hat{m} \in M, \hat{n} \in N$ such that $A = \eta^{-1} \hat{m} \hat{n} \eta t \in {\rm I} \Gamma(\P)$. Again, using $\hat{m}, \hat{n} \in K^G$ and (\ref{I-Gamma-description}) we get a contradiction for $t \neq 1$ (see \cite{P1} for details). This completes the proof of the lemma. \qed

\subsection{Integral computation}
From Lemma \ref{etam-t-in-supp-Wsharp}, we see that the integral (\ref{integral-as-a-sum-1}) is equal to
\begin{equation}\label{simplified-integral-as-a-sum-1}
Z(s) = \sum\limits_{l \geq 0} W^\#(\eta h(l,0), s) B(h(l,0)) V^{l,0}_1.
\end{equation}
Arguing as in Sect.\ 3.5 of \cite{Fu}, we get
$$V^{l,0}_1 = \frac{(1 - \Big(\frac L{\p}\Big) q^{-1}) q}{(1+q)^2(1+q^2)}q^{3l}.
$$
From Proposition \ref{values-of-B-prop} and (\ref{Wsharpformulaeq-steinberg}), we get
\begin{eqnarray}
&& B(h(l,0)) = (-\omega q^{-3})^l, \label{B-formula-for-integral}\\
&& W^\#(\eta h(l,0), s) = q^{-3(s+\frac 12)l} \omega_\tau(\varpi^{-l}) W^{(1)}(\mat{\varpi^l}{}{}{1}). \label{Wsharp-formula-for-integral}
\end{eqnarray}
Let us set $C = \frac{(1 - \big(\frac L{\p}\big) q^{-1}) q}{(1+q)^2(1+q^2)}$. We have
\begin{equation}\label{simplified-integral-as-a-sum-2}
Z(s) = C \sum\limits_{l \geq 0} (-\omega)^l q^{-3(s+\frac 12)l} \omega_\tau(\varpi^{-l}) W^{(1)}(\mat{\varpi^l}{}{}{1}).
\end{equation}
We will now substitute the value of $W^{(1)}$, from the table obtained in Sect.\ \ref{the-GL2-newform-section}, into (\ref{simplified-integral-as-a-sum-2}) for all possible $\GL_2$ representations $\tau$.
\begin{description}
\item[$\tau = \alpha \times \beta$, with $\alpha, \beta$ unramified and $\alpha \beta^{-1} \neq |\,\,|^{\pm 1}$:]  We get
\begin{equation}\label{tau-unramified-integral-formula}
Z(s) = C \frac{1}{(1+\omega \alpha(\varpi^{-1}) q^{-3s-2}) (1+\omega \beta(\varpi^{-1}) q^{-3s-2})}.
\end{equation}

\item[$\tau = \alpha \times \beta$, with $\alpha$ unramified, $\beta$ ramified and $\alpha \beta^{-1} \neq |\,\,|^{\pm 1}$:] We get
\begin{equation}\label{tau=alpha-beta-unram-ram-integral-formula}
Z(s) = C\frac{1}{1+\omega \alpha(\varpi^{-1})q^{-3s-2}}.
\end{equation}

\item[$\tau = \alpha \times \beta$, with $\alpha, \beta$ ramified and $\alpha \beta^{-1} \neq |\,\,|^{\pm 1}$:] We get
\begin{equation}\label{tau=alpha-beta-ram-integral-formula}
Z(s) = C.
\end{equation}

\item[$\tau$ supercuspidal, OR $\tau = \Omega' {\rm St}_{\GL_2}$, with $\Omega'$ ramified:] We get
\begin{equation}\label{tau=supercuspidal-OR-ramified-twisted-steinberg-integral-formula}
Z(s) = C.
\end{equation}

\item[$\tau = \Omega' {\rm St}_{\GL_2}$, with $\Omega'$ unramified:] We get
\begin{equation}\label{tau-unramified-twist-steinberg-integral-formula}
Z(s) = C \frac 1{1 + \omega \Omega'(\varpi^{-1}) q^{-3s-\frac 52}}.
\end{equation}
\end{description}
Let $\tilde\tau$ denote the contragradient of the representation $\tau$. We have the following $L$-functions for the representation $\pi = \Omega {\rm St}_{\GSp_4}$, with $\Omega$ unramified and quadratic, twisted by $\tilde\tau$.
\begin{equation}\label{twisted-L-function}
L(s, \pi \times \tilde\tau) = \left\{
                          \begin{array}{ll}
                            (1 - \Omega(\varpi) \alpha(\varpi^{-1})q^{-s-\frac 32})^{-1}(1 - \Omega(\varpi) \beta(\varpi^{-1})q^{-s-\frac 32})^{-1}, & \hbox{ if } \tau = \alpha \times \beta, \alpha, \beta \mbox{ unramified, } \\ & \qquad \alpha \beta^{-1} \neq |\,\,|^{\pm 1};\\
                            (1 - \Omega(\varpi) \alpha(\varpi^{-1})q^{-s-\frac 32})^{-1}, & \hbox{ if } \tau = \alpha \times \beta, \alpha \mbox{ unramified, } \\ 
                            & \qquad \beta \mbox{ ramified } \alpha \beta^{-1} \neq |\,\,|^{\pm 1};\\
                            (1 - \Omega(\varpi) \Omega'(\varpi^{-1})q^{-s-1})^{-1}(1 - \Omega(\varpi) \Omega'(\varpi^{-1})q^{-s-2})^{-1}, & \hbox{ if } \tau = \Omega' {\rm St}_{\GL_2}, \Omega' \mbox{ unramified; } \\
                            1, & \hbox{ otherwise. }
                          \end{array}
                        \right.
\end{equation}
From (\ref{tau-unramified-integral-formula})-(\ref{twisted-L-function}), we get the following theorem on the integral representation of $L$-functions.
\begin{theorem}\label{steinberg-int-rep-L-fn}
Let $\pi = \Omega {\rm St}_{\GSp_4}$ be the Steinberg representation of $\GSp_4(F)$ twisted by an unramified, quadratic character $\Omega$. Let $\tau$ be any irreducible, admissible representation of $\GL_2(F)$. Let $Z(s)$ be the integral defined in (\ref{ZsWdefeq}). Choose $B$ as in Sect.\ \ref{bessel-model-section} and $W^\#$ as in Sect.\ \ref{choice-of-W}. Then we have
\begin{equation}\label{int-rep-formula-in-thm}
Z(s) = Y'(s) L(3s+\frac 12, \pi \times \tilde\tau),
\end{equation}
where
$$Y'(s) = \left\{
           \begin{array}{ll}
            C (1 + \omega \Omega'(\varpi^{-1}) q^{-3s-\frac 32}), & \hbox{ if } \tau = \Omega' {\rm St}_{\GL_2}, \Omega' \mbox{ unramified; } \\
             C, & \hbox{ otherwise.}
           \end{array}
         \right.
$$
Here, $C = \frac{(1 - \big(\frac L{\p}\big) q^{-1}) q}{(1+q)^2(1+q^2)}$.
\end{theorem}

\section{Global theory}\label{global-theory-section}
In the previous section, we computed the non-archimedean integral representation of the $L$-function $L(s, \pi \times \tilde{\tau})$ for the Steinberg representation of $\GSp_4$ twisted by any $\GL_2$ representation. In \cite{Fu}, the integral has been computed for both $\pi$ and $\tau$ unramified. In \cite{PS2}, the integral has been calculated for an unramified representation $\pi$ twisted by any ramified $\GL_2$ representation $\tau$. Also, in \cite{PS2}, the archimedean integral has been computed for $\pi_\infty$ a holomorphic (or limit of holomorphic) discrete series representation with scalar minimal $K$-type, and $\tau_\infty$ any representation of $\GL_2(\R)$. In this section, we will put together all the local computations and obtain an integral representation of a global $L$-function. We will start with a Siegel cuspidal newform $F$ of weight $l$ with respect to the Borel congruence subgroup of square-free level. We will obtain an integral representation of the $L$-function of $F$ twisted by any irreducible, cuspidal, automorphic representation $\tau$ of $\GL_2(\A)$. When $\tau$ is obtained from a holomorphic cusp form of the same weight $l$ as $F$, we obtain a special value result for the $L$-function, in the spirit of Deligne's conjectures.

\subsection{Siegel modular form and Bessel model}\label{Siegel-mod-form-section}
Let $M$ be a square-free positive integer and $l$ be any positive integer. Let
$$B(M) := \{ g \in \SSp_4(\Z) : g \equiv \begin{bmatrix}\ast&0&\ast&\ast\\ \ast &\ast&\ast&\ast\\0&0&\ast&\ast\\0&0&0&\ast\end{bmatrix} \pmod{M} \}.$$
Let $F$ be a Siegel newform of weight $l$ with respect to $B(M)$. We refer the reader to Sect.\ 8 of \cite{Sa} or \cite{Sc} for definition and details on newforms with square-free level. The Fourier expansion of $F$ is given by
$$F(Z) = \sum\limits_{T>0} A(T) e^{2 \pi i \tr(TZ)},$$
where $T$ runs over all semi-integral, symmetric, positive definite $2 \times 2$ matrices. We obtain a well-defined function $\Phi= \Phi_F$ on $H(\A)$, where $\A$ is the ring of adeles of $\Q$, by
\begin{equation}\label{lift of siegel modular form to group}
 \Phi(\gamma h_{\infty} k_0) =
 \mu_2(h_{\infty})^l\det(J(h_{\infty},i1_2))^{-l} F(h_\infty \langle i1_2 \rangle),
\end{equation}
where $\gamma \in H(\Q)$, $h_{\infty} \in H^{+}(\R)$, $k_0 \in
\prod\limits_{p \nmid M}H(\Z_p) \prod\limits_{p \mid M} \I_p$. Let $V_F$ be the space generated by the right translates of $\Phi_F$ and let $\pi_F$ be one of the irreducible components. Then $\pi_F = \otimes \pi_p$, where $\pi_\infty$ is a holomorphic discrete series representation of $H(\R)$ of lowest weight $(l,l)$, for a finite prime $p \nmid M$, $\pi_p$ is an irreducible, unramified representation of $H(\Q_p)$, and for $p \mid M$, $\pi_p$ is a twist $\Omega_p {\rm St}_{\GSp_4}$ of the Steinberg representation of $H(\Q_p)$ by an unramified, quadratic character $\Omega_p$.

For a positive integer $D \equiv 0,3 \pmod{4}$, set
\begin{equation}\label{special four coeff matrix}\renewcommand{\arraystretch}{1.2}
 S(-D) = \left\{\begin{array}{l@{\qquad}l}
    \mat{\frac D4}{0}{0}{1}& \hbox{ if } D \equiv 0 \pmod{4}, \\[3ex]
 \mat{\frac{1+D}{4}}{\frac 12}{\frac 12}{1}& \hbox{ if } D \equiv 3 \pmod{4}.
 \end{array}\right.
\end{equation}
Let $L = \Q(\sqrt{-D})$ and $T(\A) \simeq
\A_L^{\times}$ be the adelic points of the group defined in
(\ref{TRthetaeq1}). Let $R(\A) = T(\A) U(\A)$ be the Bessel
subgroup of $H(\A)$. Let $\Lambda$ be a character of
\begin{equation}\label{T(A)-group}
 T(\A)/T(\Q)T(\R) \prod\limits_{p \nmid M}T(\Z_p) \prod\limits_{p \mid M} T^0_p,
\end{equation}
where, $T(\Z_p) = T(\Q_p) \cap \GL_2(\Z_p)$ and $T^0_p = T(\Z_p) \cap \Gamma^0_p$. Here $\Gamma^0_p = \{g \in \GL_2(\Z_p) : g \equiv \mat{\ast}{0}{\ast}{\ast} \pmod{p\Z_p} \}$. Note that, under the isomorphism (\ref{T-isomorphic-to-L}), $T^0_p$ corresponds to $\Z_p^\times + p \OF_{L_p}$, where $\OF_{L_p}$ is the ring of integers of the two dimensional algebra $L \otimes_{\Q} \Q_p$.
Let $\psi$ be a character of $\Q\backslash \A$ that is trivial
on $\Z_p$ for all primes $p$ and satisfies $\psi(x)=e^{-2\pi ix}$
for all $x\in\R$. We define the global
Bessel function of type $(\Lambda, \theta)$ associated to $\bar{\Phi}$ by
\begin{equation}\label{global Bessel model defn}
B_{\bar{\Phi}}(h) = \int\limits_{Z_H(\A)R(\Q)\backslash
R(\A)}(\Lambda \otimes \theta)(r)^{-1}\bar{\Phi}(rh)dr,
\end{equation}
where $\theta(\mat{1}{X}{}{1})=\psi(\tr(SX))$ and $\bar{\Phi}(h) =
\overline{\Phi(h)}$. If $B_{\bar\Phi}$ is non-zero, then $B_{\bar\phi}$ is non-zero for any $\phi \in \pi_F$. We say that $\pi_F$ has a global Bessel model of type $(\Lambda, \theta)$ if $B_{\bar\Phi} \neq 0$. We shall make the following assumption on the representation $\pi_F$.

{\bf Assumption:} $\pi_F$ has a global Bessel model
of type $(\Lambda, \theta)$ such that 
\begin{description}
\item[A1:] $-D$ is the fundamental discriminant of $\Q(\sqrt{-D})$.

\item[A2:] $\Lambda$ is a character of (\ref{T(A)-group}).

\item[A3:] For $p \mid M$, if $L \otimes \Q_p$ is split and $\Lambda_p$ is unramified, then $\Omega_p(\varpi_p)  \Lambda_p((1,\varpi_p)) \neq 1$.
\end{description}

\begin{rem}
 In \cite{Fu}, \cite{PS1}, \cite{PS2} and \cite{Sa}, non-vanishing of a suitable Fourier coefficient of $F$ is assumed, while in \cite{PS}, the existence of a suitable global Bessel model for $\pi_F$ is assumed. Let us explain the relation of the above assumption to non-vanishing of certain Fourier coefficients of $F$. Let $\{ t_j \}$ be a set of representatives for (\ref{T(A)-group}). One can take $t_j \in \GL_2(\A_{\rm f})$. Write $t_j = \gamma_j m_j \kappa_j$, with $\gamma_j \in \GL_2(\Q), m_j \in \GL_2^+(\R)$ and $\kappa_j \in \prod_{p \nmid M}\GL_2(\Z_p) \prod_{p \mid M} \Gamma^0_p$. Set $S_j := \det(\gamma_j)^{-1}\, {}^{t}\gamma_j S(-D) \gamma_j$. Note that $\{S_j\}_j$  is a \underline{subset} of the set of representatives of $\Gamma^0(M)$ equivalence classes of primitive, semi-integral positive definite $2 \times 2$ matrices of discriminant $-D$.

From \cite{Sa} or \cite{Su},
we have, for $h_{\infty} \in H^+(\R)$,
\begin{equation}\label{Bessel model arch formula}
 B_{\bar{\Phi}}(h_{\infty}) = \mu_2(h_{\infty})^l\,
 \overline{\det(J(h_{\infty},I))^{-l}}\,e^{-2\pi i\,
 \tr(S(-D)\overline{h_{\infty}\langle I\rangle})}
 \sum\limits_j \Lambda(t_j)^{-1}\overline{A(S_j)},
\end{equation}
and $B_{\bar{\Phi}}(h_{\infty}) = 0$ for $h_{\infty} \not\in H^+(\R)$. Suppose that there is a semi-integral, symmetric, positive definite $2 \times 2$ matrix $T$ satisfying
\begin{enumerate}
\item $-D = \det(2T)$ is the fundamental discriminant of $L = \Q(\sqrt{-D})$.
\item $T$ is $\Gamma^0(M)$ equivalent to one of the $S_j$.
\item The Fourier coefficient $A(T) \neq 0$.
\end{enumerate}
Then it is clear from (\ref{Bessel model arch formula}) that one can choose a $\Lambda$ such that parts ${\bf {\rm A}1, {\rm A}2}$ of the assumption are satisfied. If $M=1$ (as in \cite{Fu}, \cite{PS1}, \cite{PS2}) or, every prime $p \mid M$ is inert in $L$ (as in \cite{Sa}), then $\{S_j\}_j$ is the \underline{complete} set of representatives of $\Gamma^0(M)$ equivalence classes and hence, condition $i)$ above implies condition $ii)$ to give the assumption from \cite{Fu}, \cite{PS1}, \cite{PS2} and \cite{Sa}. We have to include part ${\bf {\rm A}3}$ of the assumption to guarantee that the Iwahori spherical vector in $\pi_p$, for $p \mid M$, is a test vector for the Bessel functional.
\end{rem}

Let us abbreviate $a(\Lambda) = \sum \Lambda(t_j) A(S_j)$. For $h \in H(\A)$, we have
$$B_{\bar\Phi}(h) = \overline{a(\Lambda)} \prod\limits_p B_p(h_p),$$
where, $B_\infty$ is as defined in \cite{PS2}, for a finite prime $p \nmid M$, $B_p$ is the spherical vector in the $(\Lambda_p, \theta_p)$-Bessel model for $\pi_p$, and for $p \mid M$, $B_p$ is the vector in the $(\Lambda_p, \theta_p)$-Bessel model for $\pi_p$ defined by Proposition \ref{values-of-B-prop} and \ref{well-defined-proposition}. For $p < \infty$, we have normalized the $B_p$ so that $B_p(1) = 1$.

\subsection{Global induced representation and global integral}\label{global-repn-integral}
Let $\tau = \otimes \tau_p$ be an irreducible, cuspidal, automorphic representation of $\GL_2(\A)$ with central character $\omega_\tau$. For every prime $p < \infty$, let $p^{n_p}$ be the conductor of $\tau_p$. For almost all $p$, we have $n_p = 0$. Set $N = \prod_p p^{n_p}$. Choose $l_1$ to be any weight occurring in $\tau_\infty$. Let $\chi_0$ be a character of $\A_L^\times$ such that $\chi_0|_{\A^\times} = \omega_\tau$ and $\chi_{0,\infty}(\zeta) = \zeta^{l_2}$ for any $\zeta \in S^1$. Here, $l_2$ depends on $l_1$ and $l$ by the formula 
$$l_2=\left\{\begin{array}{l@{\qquad\text{if }}l}
 l_1-2l&l\leq l_1,\\-l_1&l\geq l_1\end{array}\right.$$
 as in \cite{PS2}. The existence of such a character is guaranteed by Lemma 5.3.1 of \cite{PS2}. Define another character $\chi$  of $\A_L^\times$ by
$$\chi(\zeta) = \chi_0(\bar\zeta)^{-1} \Lambda(\bar\zeta)^{-1}.$$
Let $I(s,\chi_0,\chi,\tau)$ be the induced representation of $G(\A)$ obtained in an analogous way to the local situation in Sect.\ \ref{local-induced-repn-section}. We will now define a global section $f_\Lambda(g,s)$. Let us realize the representation $\tau$ as a subspace of $L^2(\GL_2(\Q)\backslash \GL_2(\A))$ and let $\hat{f}$ be the automorphic cusp form such that the space of $\tau$ is generated by the right translates of $\hat{f}$. The function $\hat{f}$ corresponds to a cuspidal Hecke newform on the complex upper half plane. Then, $\hat{f}$ is factorizable. Write $\hat{f} = \otimes \hat{f}_p$ such that $\hat{f}_\infty$ is the function of weight $l_1$ in $\tau_\infty$. For $p < \infty$, $\hat{f}_p$ is the unique newform in $\tau_p$ with $\hat{f}_p(1) = 1$. Using $\chi_0$, extend $\hat{f}$ to a function of $\GU(1,1;L)(\A)$.

For a finite prime $p$, set
$$K^G_p := \left\{
                 \begin{array}{ll}
                   G(\Z_p), & \hbox{ if } p \nmid MN; \\
                   \I \Gamma((p \OF_{L_p})^{n_{p,0}}), & \hbox{ if } p \mid M; \\
                   H(\Z_p) \Gamma((p \OF_{L_p})^{n_p}), & \hbox{ if } p \mid N, p \nmid M.
                 \end{array}
               \right.
$$
Here, in the second case, $n_{p,0} = {\rm max}(1,n_p)$.
Set $K^G(M,N) = \prod\limits_{p < \infty} K^G_p$ and let $K_\infty$ be the maximal compact subgroup of $G(\R)$.
 Let $\eta$ be the element of $G(\Q)$ defined in (\ref{eta-defn}).
Let $\eta_{M,N}$ be the element of $G(\A)$ such that the $p$-component is given by $\eta$ for $p \mid MN$ and by $1$ for $p \nmid MN$.  For $s \in \C$, define $f_\Lambda(\,\cdot\,,s)$ on $G(\A)$ by
\begin{enumerate}
 \item $f_\Lambda(g,s) = 0$ if $g \not\in M(\A)N(\A)\eta_{M,N} K_{\infty}K^G(M,N)$.
 \item If $m = m_1m_2$, $m_i \in M^{(i)}(\A)$, $n \in N(\A)$, $k =
  k_0k_{\infty}$, $k_0 \in K^G(M,N)$, $k_{\infty} \in K_{\infty}$, then
  \begin{equation}\label{section definition}
   f_\Lambda(mn\eta_{M,N} k,s) = \delta_P^{\frac 12 + s}(m)\chi(m_1)
   \hat{f}(m_2) f(k_\infty).
  \end{equation}
  Recall
  that $\delta_P(m_1m_2) =|N_{L/\Q}(m_1)\mu_1(m_2)^{-1}|^3$.
\end{enumerate}
Here, $M^{(1)}(\A)$, $M^{(2)}(\A)$, $N(\A)$ are the adelic points of the
algebraic groups defined by (\ref{M1-M2-N-defn}) and $f$ is the function on $K_\infty$ defined in \cite{PS2} as follows
$$f(g) = \left\{
         \begin{array}{ll}
           \hat{b}(g)^{l_1-l} \det(J(g,i1_2))^{-l}, & \hbox{ if } l \leq l_1; \\
           \hat{c}(g)^{l-l_1} \det(J(g,i1_2))^{-l}, & \hbox{ if } l \geq l_1.
         \end{array}
       \right.
$$
Here, we have $J(g {}^tg, i1_2) = \mat{\hat{a}(g)}{\hat{b}(g)}{\hat{c}(g)}{\hat{d}(g)}$. As in \cite{PS2}, it can be checked that $f_\Lambda$ is well-defined. For ${\rm Re}(s)$ large enough we can form the Eisenstein series
\begin{equation}\label{global-Eis-series-defn}
E(g,s;f_\Lambda) := \sum\limits_{\gamma \in P(\Q) \backslash G(\Q)} f_\Lambda(\gamma g,s).
\end{equation}
In fact, $E(g,s;f_\Lambda)$ has a meromorphic continuation to the entire plane.
In \cite{Fu}, Furusawa studied integrals of the form
\begin{equation}\label{globalintegraleq}
 Z(s,f_\Lambda,\phi)=\int\limits_{H(\Q)Z_H(\A)\backslash H(\A)}E(h,s;f_\Lambda)\phi(h)\,dh,
\end{equation}
where $\phi\in V_\pi$. Theorem (2.4) of \cite{Fu}, the ``Basic Identity'', states that
\begin{equation}\label{basicidentityeq}
 Z(s,f_\Lambda,\phi)=\int\limits_{R(\A)\backslash H(\A)}W_{f_\Lambda}(\eta h,s)B_\phi(h)\,dh,
\end{equation}
where $B_\phi$ is the Bessel function corresponding to $\phi$ and $W_{f_\Lambda}$ is the function defined by
$$ W_{f_\Lambda}(g) = \int\limits_{\Q \backslash \A}
 f_\Lambda\Big(\begin{bmatrix}1&&&\\&1&&x\\&&1&\\&&&1\end{bmatrix}g\Big)\psi(cx)dx,
 \qquad g\in G(\A).$$

The function $W_{f_\Lambda}$ is a pure tensor and we can write $W_{f_\Lambda}(g,s) = \prod_p W^\#_p(g_p,s)$. Then we see that $W^\#_\infty$ is as defined in \cite{PS2}. For a finite prime $p \nmid M$, the $W^\#_p$ is the function defined in Sect.\ 4.5 of \cite{PS2}. For $p \mid M$, the $W^\#_p$ is as in Sect.\ \ref{choice-of-W}. It follows from (\ref{basicidentityeq}) that
$$Z(s,f_\Lambda,\bar\Phi) = \prod\limits_{p \leq \infty} Z_p(s,W^\#_p, B_p),$$
where
$$Z_p(s,W^\#_p, B_p) = \int\limits_{R(\Q_p)\backslash H(\Q_p)}W^\#_p(\eta h,s)B_p(h)\,dh.$$
When $p \nmid MN, p < \infty$, the integral $Z_p$ is evaluated in \cite{Fu}. For $p = \infty$ or $p \mid N, p \nmid M$, the integral $Z_p$ is calculated in Theorems 3.5.1 and 4.4.1 of \cite{PS2}. For $p \mid M$, the integral $Z_p$ is calculated in Theorem \ref{steinberg-int-rep-L-fn}. Putting all of this together we get the following global theorem.
\begin{theorem}\label{global-thm}
 Let $F$ be a Siegel cuspidal newform of weight $l$ with respect to $B(M)$, where $l$ is any positive integer and $M$ is square-free, satisfying the assumption stated in Sect.\ \ref{Siegel-mod-form-section}. Let $\Phi$ be the adelic function
 corresponding to $F$, and let $\pi_F$ be an irreducible component
 of the cuspidal, automorphic representation generated by $\Phi$. Let $\tau$ be
 any irreducible, cuspidal, automorphic representation of $\GL_2(\A)$. Let
 the global characters $\chi$, $\chi_0$ and $\Lambda$, as well as the
 global section $f_\Lambda \in I(s,\chi,\chi_0,\tau)$, be chosen as above. Then the
 global integral (\ref{globalintegraleq}) is given by
 \begin{equation}\label{globalintreptheoremeq1}
  Z(s,f_\Lambda,\bar{\Phi})=\Big(\prod_{p\leq\infty}Y_p(s)\Big)
   \frac{L(3s+\frac12,\pi\times\tilde\tau)}
   {L(6s+1,\omega_\tau^{-1})L(3s+1,\tilde\tau\times\mathcal{AI}(\Lambda))}
 \end{equation}
 with
 \begin{equation}\label{globalintreptheoremeq2}
  Y_\infty(s)=\overline{a(\Lambda)}
   i^{l+l_2}\frac{a^+}2\pi D^{-3s-\frac l2}\,
   \frac{(4\pi)^{-3s+\frac 32 -l}}{6s+2l+l_2-1}\,\frac{\Gamma(3s+l-1+\frac{ir}2)
  \Gamma(3s+l-1-\frac{ir}2)}{\Gamma(3s+l-\frac{l_1}2 - \frac 12)}.
 \end{equation}
 Here, $\mathcal{AI}(\Lambda)$ is the automorphic representation of $\GL_2(\A)$
 obtained from $\Lambda$ via automorphic induction. The factor $Y_p(s)$ is one
 for $p \nmid MN$. For $p \nmid M, p \mid N$, the factor $Y_p(s)$ is given in Theorem 3.5.1 of \cite{PS2}. For $p\mid M$, we have $Y_p(s) = L_p(6s+1,\omega_{\tau_p}^{-1})L(3s+1,\tilde\tau_p\times\mathcal{AI}(\Lambda_p)) Y'_p(s)$, where $Y'_p(s)$ is given in Theorem \ref{steinberg-int-rep-L-fn}. The number $r$ and $a^+$ are as in the archimedean calculation in \cite{PS2}, and the constant $a(\Lambda)$ is defined in Sect.\ \ref{Siegel-mod-form-section}.
\end{theorem}

\subsection{Special values of $L$-functions}
In this section, we will use Theorem \ref{global-thm} to obtain a special value result for the $L$-function in the case that $\tau$ corresponds to a holomorphic cusp form of the same weight as $F$. Let $\Psi \in S_l(N, \chi')$, the space of  holomorphic cusp forms on the complex upper half plane $\SH_1$ of weight $l$ with respect to $\Gamma_0(N)$ and nebentypus $\chi'$. Here $N = \prod_{p} p^{n_p}$ is any positive integer and $\chi'$ is a Dirichlet character modulo $N$. $\Psi$ has a Fourier expansion
$$
 \Psi(z) = \sum\limits_{n=1}^\infty b_n e^{2 \pi i n z}.
$$
We will assume that $\Psi$ is primitive, which means that $\Psi$ is a newform, a Hecke eigenform and is normalized so that $b_1 = 1$. Let
$\omega = \otimes \omega_p$ be the character of $\A^\times / \Q^\times$ defined as the
composition
$$
 \A^\times=\Q^\times\times\R^\times_+\times\Big(\prod_{p<\infty}\Z_p^\times\Big)
 \longrightarrow\prod_{p|N}\Z_p^\times\longrightarrow\prod_{p|N}(\Z_p/p^{n_p}\Z_p)^\times
 \cong(\Z/N\Z)^\times\stackrel{\chi'}{\longrightarrow}\C^\times.
$$
 Let $K^{(0)}(N):=\prod\limits_{p|N}K^{(0)}(\p^{n_p})\prod\limits_{p\nmid N}\GL_2(\Z_p)$
with the local congruence subgroups
$K^{(0)}(\p^n)=\GL_2(\Z_p)\cap\mat{1+p^n\Z_p}{\Z_p}{p^n\Z_p}{\Z_p}$
as in (\ref{K1defeq}). Let
$K_0(N):=\prod\limits_{p|N}K_0(\p^{n_p})\prod\limits_{p\nmid N}\GL_2(\Z_p)$, where
$K_0(\p^n)=\GL_2(\Z_p)\cap\mat{\Z_p}{\Z_p}{p^n\Z_p}{\Z_p}$. Evidently,
$K^{(0)}(N)\subset K_0(N)$. Let $\lambda$ be the character of $K_0(N)$ given by
\begin{equation}\label{character-on-K0(N)-def}
 \lambda(\mat{a}{b}{c}{d}) := \prod\limits_{p | N} \omega_p(a_p).
\end{equation}
With these notations, we now define the adelic function $f_\Psi$ by
$$
 f_{\Psi}(\gamma_0 mk)=\lambda(k)\frac{\det(m)^{l/2}}{(\gamma i
 + \delta)^l}\Psi\Big(\frac{\alpha i + \beta}{\gamma i +\delta}\Big),
$$
where $\gamma_0\in\GL_2(\Q)$, $m=\mat{\alpha}{\beta}{\gamma}{\delta} \in \GL_2^{+}(\R)$
and $k\in K_0(N)$. Define a character $\chi_0$, as in the previous section, with $l_2 = -l$. Using $\chi_0$, extend $f_\Psi$ to a function on $\GU(1,1;L)(\A)$. We can take $\hat{f} = f_\Psi$ in (\ref{section definition}) and obtain the section $f_\Lambda$. Now, Lemma 5.4.2 of \cite{PS2} gives us that, for $g \in G^+(\R)$, the function $\mu_2(g)^{-l}\det(J(g,i1_2))^lE(g,s;f_\Lambda)$ only depends on $Z = g \langle i1_2 \rangle$. Let us define the function $\mathcal E$ on $\HH_2 := \{Z \in M_2(\C):\;
i(\,^{t}\!\bar{Z}-Z) \mbox{ is positive definite}\}$ by the formula
\begin{equation}\label{eis ser on siegel half plane}
 \mathcal E(Z,s) = \mu_2(g)^{-l}\det(J(g,i1_2))^l\,
 E\big(g,\frac s3 + \frac l6 - \frac 12; f_\Lambda\big),
\end{equation}
where $g \in G^+(\R)$ is such that $g\langle i1_2 \rangle = Z$. The
series that defines $\mathcal E(Z,s)$ is absolutely
convergent for ${\rm Re}(s) > 3 - l/2$ (see \cite{Kl1}). Let us assume that $l > 6$. Now, we can set $s=0$ and obtain a holomorphic Eisenstein
series $\mathcal E(Z,0)$ on $\HH_2$. Let $\Gamma^G(M,N) := G(\Q) \cap G^+(\R) K^G(M,N)$. We have $\Gamma^G(M,N) \cap H(\Q) = B(M)$. Then $\mathcal E(Z,0)$ is a modular form  of weight $l$
with respect to $\Gamma^G(M,N)$. Its restriction to $\mathfrak{h}_2$, the Siegel upper half space, is a modular form
of weight $l$ with respect to $B(M)$. By \cite{Ha}, we know that the Fourier coefficients of $\mathcal E(Z,0)$ are algebraic.

Set $V(M) := \big[\SSp_4(\Z) : B(M)\big]^{-1}$ and define, for any two Siegel modular forms $F_1, F_2$ of weight $l$ with respect to $B(M)$, the Petersson inner product by
$$\langle F_1, F_2 \rangle = \frac 12 V(M) \int\limits_{B(M)\backslash \SH_2}
  F(Z) \overline{F_2(Z)} (\det(Y))^{l-3}\,dX\,dY.
 $$
Arguing as in Lemma 5.6.2 of \cite{PS2} or Proposition 9.0.5 of \cite{Sa}, we get
\begin{equation}\label{int-as-inner-prod}
Z(\frac l6 - \frac 12, f_\Lambda, \bar{\Phi}) =  \langle \mathcal E(Z,0), F \rangle.
\end{equation}
Let $\Gamma^{(2)}(M) := \{g \in \SSp_4(\Z) : g \equiv
1 \pmod{M} \}$ be the principal congruence subgroup of $\SSp_4(\Z)$.
Let us denote the space of all Siegel cusp forms of weight $l$ with respect to
$\Gamma^{(2)}(M)$ by $S_l(\Gamma^{(2)}(M))$. For a Hecke eigenform $F \in S_l(\Gamma^{(2)}(M))$, let $\Q(F)$ be the subfield of $\C$ generated by all the Hecke eigenvalues of $F$. From \cite[p.\ 460]{Ga1}, we see that $\Q(F)$ is a totally real number field. Let $S_l(\Gamma^{(2)}(M), \Q(F))$ be the subspace of $S_l(\Gamma^{(2)}(M))$ consisting of cusp forms whose Fourier coefficients lie in $\Q(F)$. Again by \cite[p.\ 460]{Ga1}, $S_l(\Gamma^{(2)}(M))$ has an orthogonal basis $\{F_i\}$ of Hecke eigenforms $F_i \in S_l(\Gamma^{(2)}(M), \Q(F_i))$. In addition, if $F$ is a Hecke eigenform such that $F \in S_l(\Gamma^{(2)}(M), \Q(F))$, then one can take $F_1 = F$ in the above basis. Hence, let us assume that the Siegel newform $F$ of weight $l$ with respect to $B(M)$ considered in the previous section satisfies $F \in S_l(\Gamma^{(2)}(M), \Q(F))$. Then, arguing as in Lemma 5.4.3 of \cite{PS1}, we have
\begin{equation}\label{inner-prod-is-alg}
\frac{\langle \mathcal E(Z,0), F \rangle}{\langle F, F \rangle} \in \bar{\Q},
\end{equation}
where $\bar{\Q}$ is the algebraic closure of $\Q$ in $\C$. Let $\langle \Psi,\Psi \rangle_1 := (\SL_2(\Z) : \Gamma_1(N))^{-1}
\int\limits_{\Gamma_1(N)\backslash \SH_1}|\Psi(z)|^2y^{l-2}\,dx\,dy$,
where $\Gamma_1(N) := \{\mat{a}{b}{c}{d} \in \Gamma_0(N) : a, d
\equiv 1 \pmod{N} \}$. We have the following generalization of Theorem 4.8.3 of \cite{Fu}.
\begin{thm}\label{special values thm}
Let $l, M$ be positive integers such that $l > 6$ and $M$ is square-free. Let $F$ be a cuspidal Siegel newform of weight $l$ with respect to $B(M)$ such that $F \in S_l(\Gamma^{(2)}(M), \Q(F))$, satisfying the assumption from Sect.\ \ref{Siegel-mod-form-section}. Let $\Psi \in S_l(N, \chi')$ be a primitive form,
 with $N = \prod p^{n_p}$, any positive integer, and $\chi'$, any Dirichlet character modulo
 $N$.
 Let $\pi_F$ and $\tau_\Psi$ be the irreducible, cuspidal, automorphic representations
 of $\GSp_4(\A)$ and $\GL_2(\A)$ corresponding to $F$ and $\Psi$. Then
 \begin{equation}\label{special-values-eqn}
  \frac{L(\frac l2 - 1, \pi_F \times \tilde\tau_{\Psi})}{\pi^{5l-8}
  \langle F, F \rangle \langle \Psi,\Psi \rangle_1} \in \bar{\Q}.
 \end{equation}
\end{thm}
{\bf Proof.} Arguing as in the proof of Theorem 5.7.1 of \cite{PS2}, 
together with (\ref{int-as-inner-prod}) and (\ref{inner-prod-is-alg}), we get the theorem. \qed

Special value results like the one above have been obtained in \cite{BH}, \cite{Fu}, \cite{PS1}, \cite{PS2} and \cite{Sa}.


\begin{thebibliography}{99}
 \bibitem{BH} {\sc B\"ocherer, S., Heim, B.:} {\em Critical values of $L$-functions
  on $\GSp_2\times\GL_2$}. Math.\ Z. {\bf254}, 485--503 (2006)
\bibitem{Car} {\sc Cartier, P.:} {\em Representations of $p$-adic groups}. Proc. Symp. Pure Math {\bf 33}, part 1, 111--155 (1979)
 \bibitem{De1} {\sc Deligne, P.:} {\em Valeurs de fonctions $L$ et periodes d'integrales}.
  Proc. Symp. Pure Math {\bf 33}, part 2, 313--346 (1979)
 \bibitem{Fu} {\sc Furusawa, M.:} {\em On $L$-functions for $\GSp(4)\times\GL(2)$
  and their special values}. J. Reine Angew. Math. {\bf438},
  187--218 (1993)
 \bibitem{Ga1}{\sc Garrett, P.:} {\em On the arithmetic of Siegel-Hilbert cuspforms:
  Petersson inner products and Fourier coefficients}. Invent. Math. {\bf107}, no. 3,
  453--481 (1992)
 \bibitem{Ha} {\sc Harris, M.:} {\em Eisenstein series on Shimura varieties}. Ann. Math.
  {\bf 119}, 59--94 (1984)
 \bibitem{Kl1} {\sc Klingen, H.:} {\em Zum Darstellungssatz f\"ur Siegelsche Modulformen}.
  Math Z. {\bf 102}, 30--43 (1967)
 \bibitem{N-PS} {\sc Novodvorsky, M., Piatetski-Shapiro, I.:}
  {\em Generalized Bessel models for the symplectic group of rank $2$}.
  Math. USSR Sb. {\bf 19}, 243--255 (1979)
 \bibitem{P1} {\sc Pitale, A.:} {\em Steinberg representation of $\GSp_4$: Bessel models and integral representation of $L$-functions}, longer version available at www.aimath.org/$\sim$pitale/steinberg-bessel.pdf
  \bibitem{PS1} {\sc Pitale, A.; Schmidt, R.:} {\em Integral representation for
   $L$-functions for $\GSp_4\times \GL_2$}. J. Number Theory {\bf 129}, 1272--1324 (2009)
 \bibitem{PS2} {\sc Pitale, A.; Schmidt, R.:} {\em Integral representation for
   $L$-functions for $\GSp_4\times \GL_2$ II}. Preprint, 2009 available at arXiv:0908.1611
 \bibitem{PS} {\sc Pitale, A., Schmidt, R.:} {\em Bessel models for lowest weight
  representations of $\GSp(4,\R)$}. Int. Math. Res. Not. Vol. 2009, No. 7, 1159--1212 (2009)
   \bibitem{PT} {\sc Prasad, D., Takloo-Bighash, R.:} {\em Bessel models for $\GSp(4)$}.
  Preprint, 2008
\bibitem{RS}{\sc Roberts, B., Schmidt, R.:} {\em Local Newforms for $\GSp(4)$}.
  Lecture Notes in Mathematics, vol. 1918, Springer, 2007
\bibitem{Sa} {\sc Saha, A.:} {\em $L$-functions for holomorphic forms on
  $\GSp(4)\times\GL(2)$ and their special values}. Int. Math. Res. Not. Vol. 2009 No. 10, 1773--1837, 2009
   \bibitem{Sc1}{\sc Schmidt, R.:} {\em Some remarks on local newforms for $\GL(2)$}.
  J. Ramanujan Math. Soc. {\bf 17}, 115--147 (2002)
\bibitem{Sc} {\sc Schmidt, R.:} {\em Iwahori-spherical representations of $\GSp(4)$ and Siegel modular forms of degree $2$ with square-free level}.  J. Math. Soc. Japan, {\bf 57}(1), 259--293 (2005)
 \bibitem{Sh}{\sc Shimura, G.:} {\em The special values of zeta functions associated with
  cusp forms}. Comm. Pure Appl. Math. {\bf 29}, 783--804 (1976)
 \bibitem{Su} {\sc Sugano, T.:} {\em On holomorphic cusp forms on quaternion unitary
  groups of degree $2$}. J. Fac. Sci. Univ. Tokyo {\bf31},
  521--568 (1984)
\end{thebibliography}
\end{document}